\documentclass[11pt]{article}
\usepackage{amsmath}
\usepackage{lineno}
\usepackage{amssymb}
\usepackage{nicefrac} 
\usepackage{amsmath}                                                                           
\usepackage{graphicx}
\usepackage{url}
\usepackage{bm}
\usepackage[numbers]{natbib}
\numberwithin{equation}{section}
\newcommand{\bfl}{\begin{flushleft}}
\newcommand{\efl}{\end{flushleft}}

\newcommand{\be}{\begin{equation}}
\newcommand{\ee}{\end{equation}}
\newcommand{\bc}{\begin{center}}
\newcommand{\ec}{\end{center}}
\newcommand{\ben}{\begin{enumerate}}
\newcommand{\een}{\end{enumerate}}
\newcommand{\bit}{\begin{itemize}}
\newcommand{\eit}{\end{itemize}}
\newcommand{\bd}{\begin{displaymath}}
\newcommand{\ed}{\end{displaymath}}
\newcommand{\bea}{\begin{eqnarray}}
\newcommand{\eea}{\end{eqnarray}}
\newcommand{\ba}{\begin{array}}
\newcommand{\ea}{\end{array}}
\newcommand{\bfig}{\begin{figure}}
\newcommand{\efig}{\end{figure}}
\newcommand{\bpr}{\begin{proposition}}
\newcommand{\epr}{\end{proposition}}
\newcommand{\req}[1]{(\ref{eq:#1})}
\newcommand{\labeq}[1]{\label{eq:#1}}

\usepackage{amsthm}
\usepackage{ucs}
\theoremstyle{plain}
\newtheorem*{theorem*}{Theorem}

\usepackage{hyperref}
\usepackage{graphicx}
\begin{document}
\title{\sc Closed-Form Projection Method for Regularizing a Function Defined by a Discrete Set of Noisy Data and for Estimating its Derivative and Fractional Derivative}
\author{Timothy J. Burns and Bert W. Rust\\
\footnotesize Applied \& Computational Mathematics Division, NIST, Gaithersburg, Maryland, USA\\
\footnotesize \texttt{burns@nist.gov}\\ }
\date{}
\maketitle
\begin{abstract}\noindent We present a closed-form finite-dimensional projection method for regularizing a function defined by a discrete set of measurement data, which have been contaminated by random, zero mean errors, and for estimating the derivative and fractional derivative of this function by linear combinations of a few low degree trigonometric or Jacobi polynomials. Our method takes advantage of the fact that there are known infinite-dimensional singular value decompositions of the operators of integration and fractional integration. \end{abstract} 
\section{Introduction}\label{sec:IntroAbel}

In this paper, we present a finite-dimensional spectral projection method for regularizing a smooth function that is defined by a discrete set of noisy data measurements. We also present a related spectral projection method for estimating the derivative and fractional derivative of this function. Our approach takes advantage of known closed-form singular value decompositions of the integral transforms associated with integration and fractional integration.  

Our starting point is the {\em Abel transform}, an injective compact linear operator $ I_a^\mu : L^2(\Omega,w_1)\rightarrow L^2(\Omega,w_2), $
which we define by 
\begin{equation}
g(x)\ =\ I_a^\mu f(x)\ :=\  \frac{1}{\Gamma(\mu)}\int_{a}^x \left( x-y\right)^{\mu -1}f(y)dy,\label{eq:Abel}
\end{equation}
for $0<\mu<1$. Here, $\Gamma$ is the Euler gamma function, $L^2(\Omega,w_1)$ and $L^2(\Omega,w_2)$ are weighted real $L^2$ spaces on the bounded interval $\Omega = [a,b]$, with respective weight functions $w_1(x)$, $w_2(x)$, and scalar products
\be 
(f,g)_{w_1} = \int_\Omega f(x)g(x)w_1(x) dx,\hspace{5ex}(f,g)_{w_2} = \int_\Omega f(x)g(x)w_2(x) dx.\labeq{Adef}\ee
The problem of inversion of the Abel transform is a Volterra integral equation of the first kind, Abel's equation, in which the {\em data} are given by a function $g(x)$, and
the {\em source} is a function $f(x)$ to be determined.
Our goal is to find a smooth approximation of $f$, when $g$ is specified only by a finite set of measured values in $[a,b]$, which have been contaminated by random, zero-mean measurement errors, with quantified uncertainties.

Abel's equation arises in a number of important remote sensing applications, including interferometry, seismology, tomography, and astronomy; see, e.g., \citet{craigbrown}, \citet{gorenflo+vessella}, \citet{golberg}, \citet{bracewell}. The operator $I_a^\mu$ of \req{Abel} can be interpreted as a {\em fractional integral operator}, for $\mu\in(0,1)$ \citep[Sec.\ 1.1]{gorenflo+vessella}. 
For $\mu\in(0,1)$, it can be shown that \citep[Sec.\ 1.2, 1.A]{gorenflo+vessella}, if $g(x)$ is absolutely continuous on $[a,b]$,
then \req{Abel} has a unique solution $f(x)\in L^1(\Omega)$, which is given by
\begin{eqnarray}
f(x) &=& D_a^\mu g(x) := \frac{1}{\Gamma({1-\mu})}\frac{d}{dx}\int_a^x  \left( x-y\right)^{-\mu}g(y)dy \label{eq:invAbel1}\\
&= & \frac{d}{dx} I_a^{1-\mu}g(x) = (I_a^\mu)^{-1} g(x).\nonumber 
\end{eqnarray}
Thus, the inverse of $I_a^\mu$ can be interpreted as a {\em fractional derivative operator} for $\mu\in(0,1)$  \citep[Sec.\ 1.2]{gorenflo+vessella},
\begin{equation}
 D_a^\mu g(x) = \frac{d}{dx}I_a^{1-\mu}g(x).\labeq{invAbel2}\ee
If $g(a)=0$, and in addition $g(x)\in C^1[a,b]$, then by partial integration,
\begin{equation}
f(x) = D_a^\mu g(x) = \frac{1}{\Gamma({1-\mu})}\int_a^x  \left( x-y\right)^{-\mu}\frac{dg}{dy}(y)dy \label{eq:invAbel3}.
\end{equation}

The presence of the differentiation operator in \req{invAbel1} and \req{invAbel3} suggests that there will be difficulties similar to those encountered when trying to estimate the derivative of a function that is specified by noisy data, i.e., in trying to solve
\be g(x) =  I_a^1 f(x) := \int_a^x f(x) dx \labeq{integration}\ee
for $f(x)$ numerically, when $g(x)$ is given as a discrete set of data with measurement uncertainties. Like the integration operator $I_a^1$ \citep[Thm.\ 2.28]{kress}, the Abel transform $I_a^\mu$ is a compact linear operator acting on and into infinite-dimensional Hilbert spaces \citep[Thm.\ 4.3.3]{gorenflo+vessella}. Therefore, its inverse $D_a^\mu = (I_a^\mu)^{-1}$ cannot be continuous \citep[Sec.\ 2.5]{kress}. Thus, the problem of 
inverting the Abel transform is also ill-posed. 

Without loss of generality, assume for the moment that $a=0$ in \req{Abel}, and denote its convolution kernel function by $k(x-y)$, i.e.,
\be 
k(x-y) := \frac{1}{\Gamma(\mu)}(x-y)^{\mu-1}. \labeq{kernel}\ee 
Following \citet[Sec.\ 4.4]{craigbrown}, we take the Laplace transform of both sides in \req{Abel}, and get that 
\begin{equation}
\bar{g}(s) = \bar{k}(s) \bar{f}(s),\labeq{laptransfAbel}\ee
where
\begin{equation}
\bar{k}(s) = {1}/{{s}^{\mu}}.\labeq{laptranskernel}\ee
Let $s=\gamma +i\omega$, where $\gamma>0$ is a fixed real number. It follows from \req{laptranskernel} that, for large 
$|\omega|$, the point spread or filtering property of the kernel  behaves like
\be |\bar{k}(s)|\ \sim {1}/{{|\omega|}^\mu}.\labeq{filter}\ee 
Thus, the closer $\mu$ is to $1$, the more the kernel in \req{Abel} smoothes high-frequency content in the source function. It follows that if the data function is contaminated by high-frequency noise, it will not lie in the range of the operator $ I_a^\mu$, so that the Abel equation becomes more ill-posed as $\mu\rightarrow 1$, and the most unstable case is $\mu=1$, i.e., numerical differentiation of noisy data. 

Consider the following example of \req{integration} (\citet[Sec.\ 1.3]{craigbrown}),
\bea g(x) &=&  1-\exp(-\alpha x)+\beta\sin(\omega x),\labeq{gcraig}\\[5pt]
     f(x) &=&  \alpha\exp(-\alpha x)+\beta\omega\cos(\omega x),\labeq{fcraig}\eea
for $0\le x\le 2$. 
Denote by $g_a$ and $f_a$ the respective functions \req{gcraig} and \req{fcraig} evaluated with the parameters $\alpha = 0.8$, $\omega = 20$, and $\beta = 0.04$, and by $g_b$ and $f_b$ the respective functions evaluated with the parameters $\alpha = 0.8$, $\omega = 20$, and $\beta = 0$; see Fig.~\ref{fig:CraigExact0_2}. Note that, whereas the two functions $g_a$ and $g_b$ differ only by a small-amplitude oscillation of relatively low frequency, corresponding to a maximum pointwise difference of about $4\%$, the corresponding source functions $f_a$ and $f_b$ differ pointwise by a maximum of more than $80\%$. Furthermore, whereas the maximum difference between the data functions remains the same as $\omega$ increases, this is not the case for the source functions, where the maximum separation between $f_a$ and $f_b$ increases linearly with $\omega$. 
It follows that, for this example, while the forward problem of model fitting, in which source functions are compared as to how well they can predict a given set of measurement data, is insensitive to small high-frequency perturbations in $f$, the inverse problem, of estimating the source function given a set of data measurements, is extremely sensitive to the presence of high-frequency content in the data, even if it is of relatively small amplitude. That this is generally the case for \req{Abel} and \req{integration} follows from the singular value decompositions of the associated integral operators; see Sec.~\ref{sec:projection}.
Furthermore, by \req{laptranskernel}, small-amplitude high-frequency content in a set of experimental measurements of the function $g$, that is based on either of the transforms \req{Abel} or \req{integration}, for some source $f$, would be expected to correspond to noise, rather than signal, in the data.

\begin{figure}[tbp] 
  \centering
  \includegraphics[width=5in,height=2.82in,keepaspectratio]{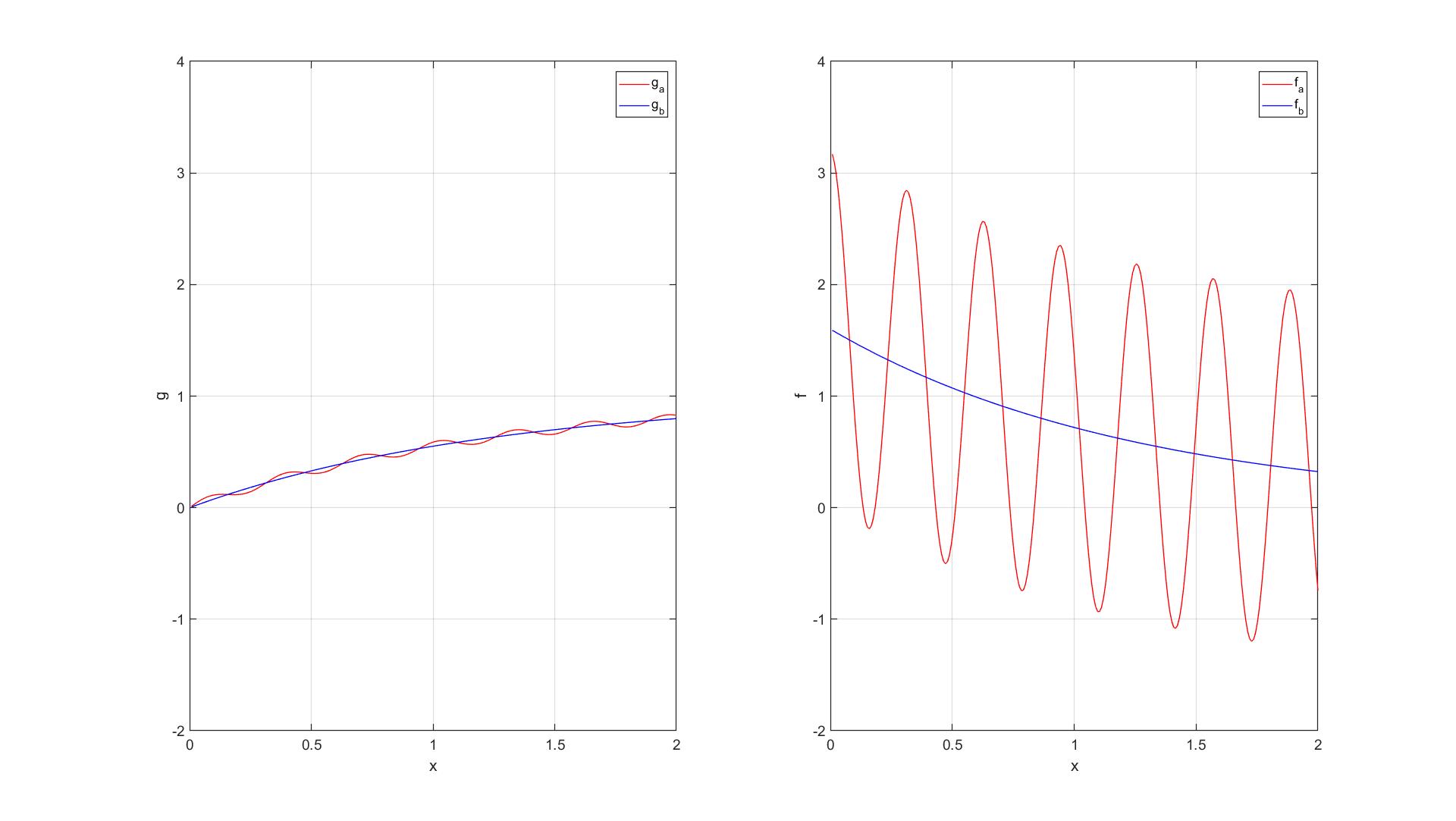}
  \caption{Plots of example of Craig and Brown \req{gcraig} and \req{fcraig}. Subscripts $a$ and $b$ refer to $\beta = 0.04$ and $\beta = 0$, respectively.}
  \label{fig:CraigExact0_2}
\end{figure}

As will be discussed in Sec.~\ref{sec:rustreg} and \ref{sec:newreg}, our approach to regularizing the data function $g$, i.e., for separating signal from noise in the data, prior to obtaining an estimate of the source function $f$, is based on a global approximation method. The important point that we want to make here is that our approach to smoothing the data is essentially independent of any specific integral equation. 
The regularization method that we will present is an extension of the work of  \citet{BertIR}
on the highly ill-posed problem
of estimating the source function in Fredholm integral equations of the first kind, 
\be  g(x) = Af(x) := \int_a^b K(x,y)f(y)dy,\labeq{fredholm}\ee
with a smooth kernel $K(x,y)$ (see \citet{craigbrown}, \citet{wing}, \citet{kress}).
Rust developed a global projection method, which is based on the numerical singular value decomposition (SVD) of a discrete version of the operator $A$ in \req{fredholm}, for the regularization of data, prior to 
computing an estimate of $f(x)$ in \req{fredholm}, using a statistical linear regression model, which treats the residual as a realization of a white noise time series.  

Rust's method will be discussed in more detail in the next section. In Sec.~\ref{sec:projection}, we will take advantage of the fact that there are known singular value decompositions for the operators of integration and fractional integration in an infinite-dimensional Hilbert space setting. By approximating $g$ in a finite-dimensional subspace of the range of the integral operator, we will show how to obtain a closed-form finite-dimensional estimate for the source function $f$ in a finite-dimensional subspace of the domain of the operator. We will then discuss our approach to obtaining a regularized approximation of $g$, prior to estimating $f$, in the context of the infinite-dimensional SVD's, in Sec.~\ref{sec:newreg}. At the end of Sec.~\ref{sec:newreg}, we will show how this regularization method can be ``decoupled" from the SVD, which will provide a method for separating higher-frequency noise from low-frequency signal in measurement data from more general applications, under the assumption that the process that produces the data is smoothing, so that it preferentially damps higher-frequency content in the data. In Sec.~\ref{sec:numerical}, we will present some computational results that show the usefulness of our regularization method. Some discussion and concluding remarks will be given in Sec.~\ref{sec:conclusions}.

\section{Rust's Regularization Method}\label{sec:rustreg}

As discussed in the Introduction, Rust developed a method for separating signal from noise in noisy data, prior to obtaining an estimate of the source function in an overdetermined linear regression model of a Fredholm integral equation of the first kind. Subsequently, Rust's method has influenced the work of others on first-kind Fredholm integral equations; see, e.g., \citep{ha+ki+kj}. We will show in what follows that this approach to regularization has much wider application to ill-posed problems. 

In this section, we will outline Rust's method by finding an approximate solution of the simple-looking Volterra integral equation of the first kind (see \req{integration}),  
\be g(x) \ =\  \int_a^x  f(y)\,dy,\hspace{5ex}-\infty <a\le x \le b<\infty.\labeq{integration2}\ee
In \req{integration2}, data represented by the function $g$ evaluated at $x$ are the result of equal contributions from the source function $f$, but only over the subset $a\le y\le x$ of its domain of definition. This is in contrast to a Fredholm integral equation of the first kind \req{fredholm}, in which the source $f$ contributes to the data 
$g$ over the entire interval on which the source is defined.
 
We assume that the integral equation \req{integration2} has been nondimensionalized. We also assume that $g(a)=0$, and that the function $g$ is defined by a discrete set of measured values, $\mathbf{g}:=(g_1,\ldots,g_m)^T$, corresponding to the mesh, $a<x_1,\ldots,x_m =b$, which have been contaminated by random, zero-mean errors 
\be \boldsymbol{\epsilon}:= (\epsilon_1,\ldots,\epsilon_m)^T,\labeq{errors}\ee 
such that 
\begin{equation} 
\mathbf{\mathcal{E}} \left(\boldsymbol{\epsilon}\right)\ =\ \mathbf{0}, \hspace{5ex}\mathbf{\mathcal{E}}\left(\boldsymbol{\epsilon \epsilon}^T \right)\ =\ \mathbf{S}^2.\labeq{calE}
\end{equation}
Here, $\mathbf{\mathcal{E}}$ is the expectation operator, $\mathbf{0}$ is the $m$-dimensional zero vector, and $\mathbf{S}^2$ is the $m\times m$ positive definite variance-covariance matrix for $\boldsymbol{\epsilon}$.  It is also assumed that the measurement errors are statistically independent, so that
\begin{equation} 
\mathbf{S}^2 \ =\  \mathbf{diag}(s_1^2,\ldots,s_m^2),\labeq{S2}\ee
where $s_1,\,\ldots,\, s_m$ are the estimated standard deviations of the errors. In experimental science, these are an integral part of the measurement process. They are typically reported as $\pm1\,s_j$ or $\pm2\,s_j$ error bars on published plots of data \citep{gum}. The $x_j$ need not be equally spaced. We take the point of view that we have no control over the mesh size $m$ or the mesh spacing, so that, for example, we may not assume that the $x_j$ are Chebyshev nodes. For simplicity, we will assume that the mesh points $x_j$ are equally spaced. 

\subsection{Discrete Problem}\label{subsec:discrete}
The inverse problem can then be written as a sequence of integral equations,
\be g_j \ =\  \int_a^{x_j}  f(y)\,dy + \epsilon_j,\hspace{5ex} j=1,\ldots,m.\labeq{g_j}\ee
We fully discretize the problem by 
approximating the integrals \req{g_j} using the midpoint rule, and get
\be g_j \ =\   h\sum_{i=1}^j \tilde{f}_{i-\nicefrac12} + \epsilon_j ,\hspace{3ex} j=1,\ldots,m,\labeq{midpt}\ee 
where  $h=(b-a)/m$, $x_{j-\nicefrac12} = (x_{j-1}+x_j)/2,$ and $\tilde{f}_{i-\nicefrac12}=f(x_{i-\nicefrac12}).$  
We assume that $h$ is small enough to capture the details of $g$, and that 
\be s _j \gg h^2,\hspace{5ex}j=1,\ldots,m.\labeq{epsh2}\ee 
The requirement \req{epsh2} guarantees that the discretization error of the midpoint rule, which is O$(h^2)$, is much smaller than the measurement errors.
Henceforth, we will neglect the discretization error.

Problem \req{midpt} can be written as
\be {\mathbf{g}} \ =\  \mathbf{L}\tilde{\mathbf{f}}+\boldsymbol{\epsilon},\hspace{5ex}\boldsymbol{\epsilon}\sim \mbox{N}(\boldsymbol{0},\mathbf{S}^2),\labeq{Lfeps}\ee
where 
$\mathbf{L}$ is the $m\times m$ lower-triangular, non-singular matrix
\be
\mathbf{L} \ =\  
h\begin{pmatrix}
    1 & &  & & \\
    1 & 1 &  & & \\
    1 & 1 & 1 & & \\
    \vdots & \vdots&\vdots & \ddots   \\
    1 & 1 & 1   &  \cdots  & 1 \\
\end{pmatrix}.\labeq{Lmatrix}\ee
In the present case, there is no overdetermined linear regression problem to solve. It is straightforward to verify that the inverse of $\mathbf{L}$ is given by
\be
\mathbf{L}^{-1} = 
\frac{1}{h}\begin{pmatrix}
    1 & &  &  \\
    -1 & 1 &  &  \\
     &  -1 & 1   &  \\
   &  &  \ddots &     \ddots  & \\
   &  &  &            -1  &    1    \\
\end{pmatrix}
,\labeq{Linvmatrix}\ee
so that an estimate of the source function is obtained by first-order finite-dif\-fer\-enc\-ing of the noisy data.

\subsection{Need for Regularization}\label{subsec:need}
For Case $(a)$ of \req{gcraig}-\req{fcraig} (see Fig.~\ref{fig:CraigExact0_2}), with $m=250$, $0\le x\le2$, $\alpha = 0.8$, $\omega = 20$, and $\beta = 0.04$, we perturb each data point $g(x_j)$ as follows. Obtain a pseudorandom sample from a standard normal distribution, using the Matlab function $\mathbf{randn}$, multiply it by a fixed standard deviation of $s = 0.05$, and then add this number to $g(x_j),\ j=1,\ldots,m$. The resulting estimate $\hat{\mathbf{f}}$ of $f$ is plotted on the right-hand side of Fig.~\ref{fig:gpertfpert}.
The estimate of the source function is dominated by amplified noise in the data, which is a characteristic feature of an ill-posed problem. It is apparent that some kind of regularization of $\mathbf{g}$, i.e., separation of signal from noise in the data, is necessary prior to obtaining an estimate of $f$.

\begin{figure}[tbp] 
  \centering
  \includegraphics[width=5in,height=2.35in,keepaspectratio]{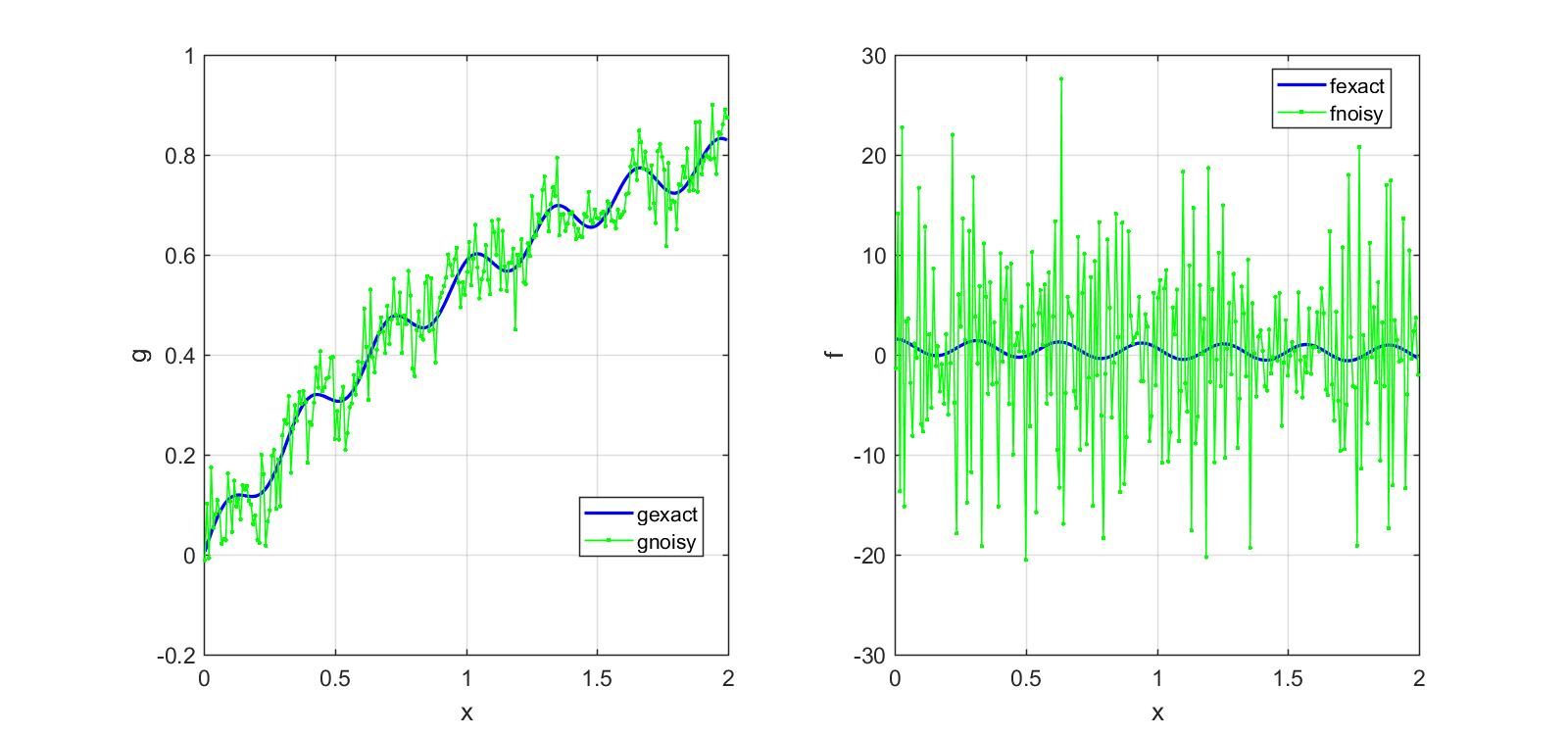}
  \caption{Left: Exact and perturbed data function $g$ \req{gcraig}, with parameters $\alpha=0.8,\ \beta=0.04,\ \omega=20,$ and $m=250$, $s_j =0.05,\ j=1\ldots,m$; noise generated using Matlab function $\mathbf{randn}$.
Right: Exact and noisy source function $f$ \req{fcraig}; $\mathbf{fnoisy} = \mathbf{L}^{-1}\mathbf{gnoisy}$.}
  \label{fig:gpertfpert}
\end{figure}

\subsection{Scaling the Problem}\label{subsec:scaling}
Rust's method for regularizing the data begins with scaling the problem. Multiplying both sides of \req{Lfeps} by 
$\mathbf{S}^{-1}$, we get
\be \mathbf{b}\ =\ \mathbf{M}\tilde{\mathbf{f}}+\boldsymbol{\eta},\labeq{disc2}\ee
where
\be \mathbf{b}\ =\ \mathbf{S}^{-1}\mathbf{g},\hspace{5ex}\mathbf{M}\ =\ \mathbf{S}^{-1}\mathbf{L},\hspace{5ex}
\boldsymbol{\eta}\ =\ \mathbf{S}^{-1}\boldsymbol{\epsilon}.\labeq{scaling}\ee
For the scaled problem \req{disc2}, it follows by a standard theorem of multivariate statistics \citep[Thm 2.4.4]{anderson}, that
\be \boldsymbol{\eta}\sim \mbox{N}\left(\mathbf{0},\mathbf{I}_m\right),\labeq{etadist}\ee
or
\be \mathbf{b}\sim \mbox{N}\left(\mathbf{M}\tilde{\mathbf{f}},\mathbf{I}_m\right).\labeq{bdist}\ee

The considerable advantage that is gained by this scaling can be seen as follows. 
Let $\hat{\mathbf{f}}$ be a discrete estimate of the unknown source function $\tilde{\mathbf{f}}$, and
let
\be \hat{\mathbf{r}} \ :=\  {\mathbf{b}}-\mathbf{M}\hat{\mathbf{f}}\labeq{hatr}\ee
be the corresponding residual vector. Since the regression model can also be written
\be  \boldsymbol{\eta} \ =\  \mathbf{b}-\mathbf{M}\tilde{\mathbf{f}},\labeq{etaeq}\ee
it is clear that {\em an estimate} $\hat{\mathbf{f}}$ {\em is acceptable only if} $\hat{\mathbf{r}}$ {\em is a plausible sample from the} 
$\boldsymbol{\eta}$ {\em distribution.} An important feature of this approach is that the elements of the residual vector $\hat{\mathbf{r}} $ and 
the true residual vector $\boldsymbol{\eta}$ are considered to be discrete time series, in which the component number is treated as the time variable.

\subsection{Diagnostics Based on Statistical Properties of the Residual}\label{subsec:diagnostics}

When the source $f$ is unknown, the residual vector provides the only objective guide for assessing the quality of an estimate of this function. In the numerical example in Sec.\ \ref{subsec:need} (see Fig.~\ref{fig:gpertfpert}), $\|\hat{\mathbf{r}}\|^2_2 \approx 0$, which is much too small. Rust \citep{BertIR} (see also \citet{rust2000},\  \citet{rust+oleary}) proposed several diagnostics for judging the acceptability of an estimate $\hat{\mathbf{f}}$ with residual $\hat{\mathbf{r}}$. 
Since the $\eta_j$ are statistically independent and identically normally distributed, the residuals for an acceptable estimate of the source function should constitute a realization of such a time series. 

The first diagnostic is a version of {\em Morozov's Discrepancy Principle} \citep{morozov1966}, which states that the size of the residual should be comparable to the magnitude of the error in the data. Rust formulated this principle, in the present context, as follows.
Since $\mathbf{b}-\mathbf{M}\tilde{\mathbf{f}}\sim \mbox{N}\left(\mathbf{0},\mathbf{I}_m\right),$ it follows from a standard statistical theorem \citep[p 140]{hogg+craig} that
\be \|{\mathbf{b}-\mathbf{M}\tilde{\mathbf{f}}}\|_2^2 \ \sim \chi^2(m),\labeq{chi2}\ee
where $\chi^2(m)$ denotes the chi-square distribution with $m$ degrees of freedom. It follows that
\be \mathbf{\mathcal{E}} \{ \| \mathbf{b}-\mathbf{M}\tilde{ \mathbf{f} } \|_2^2 \} \ =\  m,\hspace{5ex}
\mbox{Var} \{ \| \mathbf{b}-\mathbf{M}\tilde{ \mathbf{f} } \|_2^2 \} \ =\  2m.\labeq{ssrbounds}\ee
These two quantities provide rough bounds for the sum of squared residuals that might be expected from a reasonable estimate
$\hat{ \mathbf{f} }$ of $\tilde{ \mathbf{f} }$. An estimate such that
$ m-\kappa\sqrt{2m}\le\| \mathbf{b}-\mathbf{M}\hat{ \mathbf{f} } \|_2^2 \le m+\kappa\sqrt{2m}$, 
with $\kappa=1$,
would be quite reasonable, but any $\hat{ \mathbf{f} }$ whose sum of squared residuals falls outside the bounds with $\kappa=2$
would be suspect. This leads to 

\vspace{1em}
\noindent  
\textbf{Diagnostic 1.}
\be    m-2\sqrt{2m}\ \le\ \| \hat{ \mathbf{r} } \|_2^2\ \le\ m+2\sqrt{2m}.\labeq{roughbounds}\ee

\vspace{1em}
\noindent \textbf{Diagnostic 2}. {\em The elements of $\boldsymbol{\eta}$ are drawn from a} $\mbox{N}(0,1)$ {\em distribution; therefore, a graph of the
elements of $\hat{\mathbf{r}}$ should look like samples from this distribution. In fact, a histogram of the
entries of $\hat{\mathbf{r}}$ should look like a bell curve. } 

\vspace{1em}
A quantititative method for evaluating Diagnostic $2$ is the chi-square goodness-of-fit test for a normal distribution \citep{sn+co}. Software for evaluating this criterion is available in the Matlab Statistics toolbox (Matlab function $\mathbf{chi2gof}$).

\vspace{1em}
The third diagnostic is based on the periodogram of the residuals, which is an estimate of the spectral density of the residual time series, i.e., it is an estimate of how the total variance in the series is distributed in frequency. A more detailed discussion of this diagnostic is given by Fuller \citep[Ch.\ 7]{fuller}. 

To define the periodogram, we first need to define the discrete Fourier transform (DFT). The {\em discrete Fourier transform} of the time series data set $(r_1,\ldots,r_m)$ is defined to be the set of complex numbers $(\mathcal{R}_{0},\ldots,\mathcal{R}_{m-1})$, where
\bd \mathcal{R}_j = \sum_{t=1}^m r_t\, \exp(-i\, 2\pi \nu_j t),\hspace{3ex}\nu_j= j/m,\hspace{3ex}j=0,\ldots,m-1,\ed
where $i=\sqrt{-1}$; the {\em Fourier frequencies} $\nu_j = j/m$, $0\le \nu_j \le 1$, are harmonics of the fundamental frequency $1/m$.

The {\em periodogram} is defined by
\be \mathcal{P}_{j}= \frac{1}{m}  |\mathcal{R}_{j}|^2,\hspace{3ex}j=0,\ldots,q,\labeq{periodog}\ee 
where $q$ is the smallest integer greater than or equal to $(m-1)/2$.
The presence of a sinusoidal cycle in $(r_1,\ldots,r_m)$ is indicated by a peak in a plot of $\mathcal{P}_j$ against the frequencies $\nu_j = j/m$, $j=0,\ldots,q$, $0\le \nu_j \le 0.5$. The amplitude of the peak gives an estimate of the power in that cycle. 

The {\em cumulative periodogram} is defined by $\mathcal{C}_0,\mathcal{C}_1,\ldots,\mathcal{C}_{q}$, where $\mathcal{C}_0  = 0$, and
\be \mathcal{C}_j =  \left[ \sum_{k=1}^q \mathcal{P}_k \right]^{-1} \sum_{k=1}^{j} \mathcal{P}_k ,\ \ j=1,\ldots,q.\labeq{cumper}\ee
We note that the $\mathcal{C}_j$ are nondecreasing as $j$ increases, and that $\mathcal{C}_{q} =1$.
It can be shown that, for a realization of a white noise time series, a plot of the cumulative periodogram against the frequencies $\nu_j = j/m,\ j=0,\ldots,q$, will fluctuate about the straight line $y=2\nu$, from $[0,0]$ to $[0.5,1]$, which corresponds to an ideal white noise series in which the variance is uniformly distributed. 
This line has length 
\be\mathcal{L} = \sqrt{(0.5)^2+(1.0)^2} \approx 1.1180.\labeq{calL}\ee
A useful check is to compare $\mathcal{L}$ with 
\be \mathcal{L_C} = \sum_{k=1}^q \sqrt{(\mathcal{C}_{k}-\mathcal{C}_{k-1})^2+(\nu_k - \nu_{k-1})^2}.\labeq{LC}\ee
The null hypothesis that the residual series $(r_1,\ldots,r_m$) represents a realization of white noise can be tested by plotting the two straight lines $y=2\nu \pm \delta,$ where $\delta$ is the $5\%$ point for the Kolmogorov-Smirnov statistic for a sample of size $m-1$. These two lines define a $95\%$ confidence band for white noise. 

\vspace{1em}
\noindent \textbf{Diagnostic 3}. {\em No more than $5\%$ of the ordinates of a plot of the cumulative periodogram should lie outside the $95\%$ confidence band.}

\vspace{1em}
It is convenient to use the FFT algorithm to compute the DFT and the periodogram. We do this by first padding the residual series by zeros, i.e., by choosing $M$ to be a convenient power of two, such that $M \ge m$, then defining $r_j = 0$, for $j = {m+1},\ldots,M $, and then re-defining the DFT by
\be  \mathcal{R}_j = \sum_{t=1}^M r_t\, \exp(-i\, 2\pi \nu_j t),\hspace{3ex}\nu_j= j/M,\hspace{3ex}j=0,\ldots,{M}/{2},\labeq{DFT}\ee
so that the new Fourier frequencies are harmonics of the fundamental frequency $1/M$. The zero-padding increases the density of the frequency mesh, but it does not change the value of the Fourier transform at any given frequency. The periodogram and the cumulative periodogram are computed as before, except that now $q= M/2$.

\subsection{Separating Signal from Noise in the Data}\label{subsec:rustmethod}
In the setting of overdetermined least-squares regression, Rust regularized the data by means of a projection method. For the case of the re-scaled problem \req{disc2} that we consider here, where $\mathbf{M}$ is an invertible $m\times m$ square matrix, the method is as follows. First, compute the singular value decomposition of $\mathbf{M}$,
\be \mathbf{M}=\mathbf{U} \boldsymbol{\Sigma}\mathbf{V}^T.\labeq{finSVD}\ee
Here, each of the three matrices on the right-hand side is $m\times m$. $\mathbf{U}$ is orthogonal, and its columns $\{\boldsymbol{u}_1,$ $\ldots,$ $\boldsymbol{u}_m\}$ span the range of $\mathbf{M}$. In the examples that we have examined, if one plots the corresponding $\{\mathbf{u}_j,\ j=1,\ldots,m\}$, the vectors appear to become more and more oscillatory as the index $j$ increases. In Sec.~\ref{sec:newreg}, 
we will show that this behavior is to be expected for the operators of integration \req{integration} and fractional integration \req{Abel}, in an infinite-dimensional Hilbert space setting.

The central matrix in \req{finSVD} is 
diagonal, 
$\boldsymbol{\Sigma}=
\mbox{diag}(\sigma_1, \ldots,\sigma_m)$. The entries along the main diagonal are the singular values of       
$\mathbf{M}$, all positive numbers, arranged in descending order. $\mathbf{V}$ is also an orthogonal matrix, with columns $\{\boldsymbol{v}_1,\ldots,\boldsymbol{v}_m\}$ that span the domain of $\mathbf{M}$. 

The next step is to project the data vector $\mathbf{b}$ onto each of the orthonormal basis vectors $\boldsymbol{u}_j$ of the range of $\mathbf{M}$, by multiplying both sides of \req{disc2} by the matrix $\mathbf{U}^T$, so that
\be \mathbf{U}^T\mathbf{b} = \boldsymbol{\Sigma}\mathbf{V}^T \tilde{\mathbf{f}}+\mathbf{U}^T\boldsymbol{\eta}.\labeq{UTL}\ee
Because $\mathbf{U}$ is an orthogonal matrix, the distribution of each component $\left(\mathbf{U}^T\boldsymbol{\eta}\right)_j$ of the $\mathbf{U}^T\boldsymbol{\eta}$ vector is also $\sim \mbox{N}(0,1)$, and  
\be \| \mathbf{b}-\mathbf{M}\tilde{ \mathbf{f} } \|_2^2 = \|\mathbf{U}^T\mathbf{b}- \boldsymbol{\Sigma}\mathbf{V}^T \tilde{\mathbf{f}}\|_2^2  \sim \chi^2(m).\labeq{UTchi2}\ee

This suggests a practical method for separating signal from noise in the data vector $\mathbf{b}$. 
Project the data vector onto the $m$ column vectors of $\mathbf{U}$, 
\be \left(\mathbf{U}^T\mathbf{b}\right)_j = \sum_{k=1}^m b_k \boldsymbol{u}_j^T \mathbf{i}_k,\hspace{3ex}j=1,\ldots,m,
\labeq{bproject}\ee 
where $\{\mathbf{i}_1,\ldots,\mathbf{i}_m\}$ is the standard basis of $\mathbb{R}^m$.
To simplify the notation, let 
\be \mathbf{a} := \mathbf{U}^T\mathbf{b}.\labeq{adef}\ee

\vspace{1ex}\noindent
{\em Since most experimentalists would be reluctant to claim that a measured value is significantly different from zero if its magnitude does not exceed three standard deviations for the estimated error in the measurement, separate the projected data  
into signal and noise, respectively, $\mathbf{a}=\mathbf{a}_S + \mathbf{a}_N$, according to whether or not 
$|{a}_j| > \tau$, for $j=1,\ldots,m$, where $\tau = 3$.} 

\vspace{1ex}\noindent
This is a rule of thumb that we have found to work well in the examples that we have studied.\ {\em Start with $\tau =3$, and then, if necessary, adjust the value of $\tau$ upward or downward, to ensure that Morosov's discrepancy principle \req{roughbounds} is satisfied}. 
We then have the truncation of the scaled data vector,
\be\mathbf{b}\ =\ \mathbf{U}\mathbf{a} \ =\  \mathbf{U}\mathbf{a}_S + \mathbf{U}\mathbf{a}_N \ =\ \mathbf{b}_S+ \mathbf{b}_N,\labeq{rustbtrunc}\ee 
and of the original data vector,
\be\mathbf{g}\ =\ \mathbf{S}\mathbf{b} \ =\  \mathbf{S}\mathbf{b}_S + \mathbf{S}\mathbf{b}_N \ =\ \mathbf{g}_S+ \mathbf{g}_N.\labeq{rustgtrunc}\ee

For the problem of Craig and Brown in Fig.~\ref{fig:gpertfpert}, with $\tau = 3$ as the boundary between signal and noise in $\mathbf{a}$, the lower and upper $\chi^2$ bounds \req{roughbounds} are given by $205.3$ and  $294.7$, respectively. We find that $\mathbf{a}_S$ consists of six components, corresponding to the indices $1,\     2,\     3,\    13,\    24,$ and $192$, with the respective values of $-174.4,\ 
-3.5,\ -4.2,\ -8.1,\ 3.1,$ and $3.6$; see Fig.~\ref{fig:aprojMidPt}. The sum of squared residuals $\|\mathbf{a}_N\|_2^2 = 258.1$, which satisfies \req{roughbounds}. It turns out that all of the column vectors of $\mathbf{U}$ that correspond to the six signal components of $\mathbf{a}$ are oscillatory. The first five of these oscillate with increasing frequency; see Fig.~\ref{fig:aprojMidPt}. However, the column of $\mathbf{U}$ that corresponds to component number $192$ oscillates with a much higher frequency than the first five.

\begin{figure}[tbp] 
  \centering
  \includegraphics[width=5in,height=2.52in,keepaspectratio]{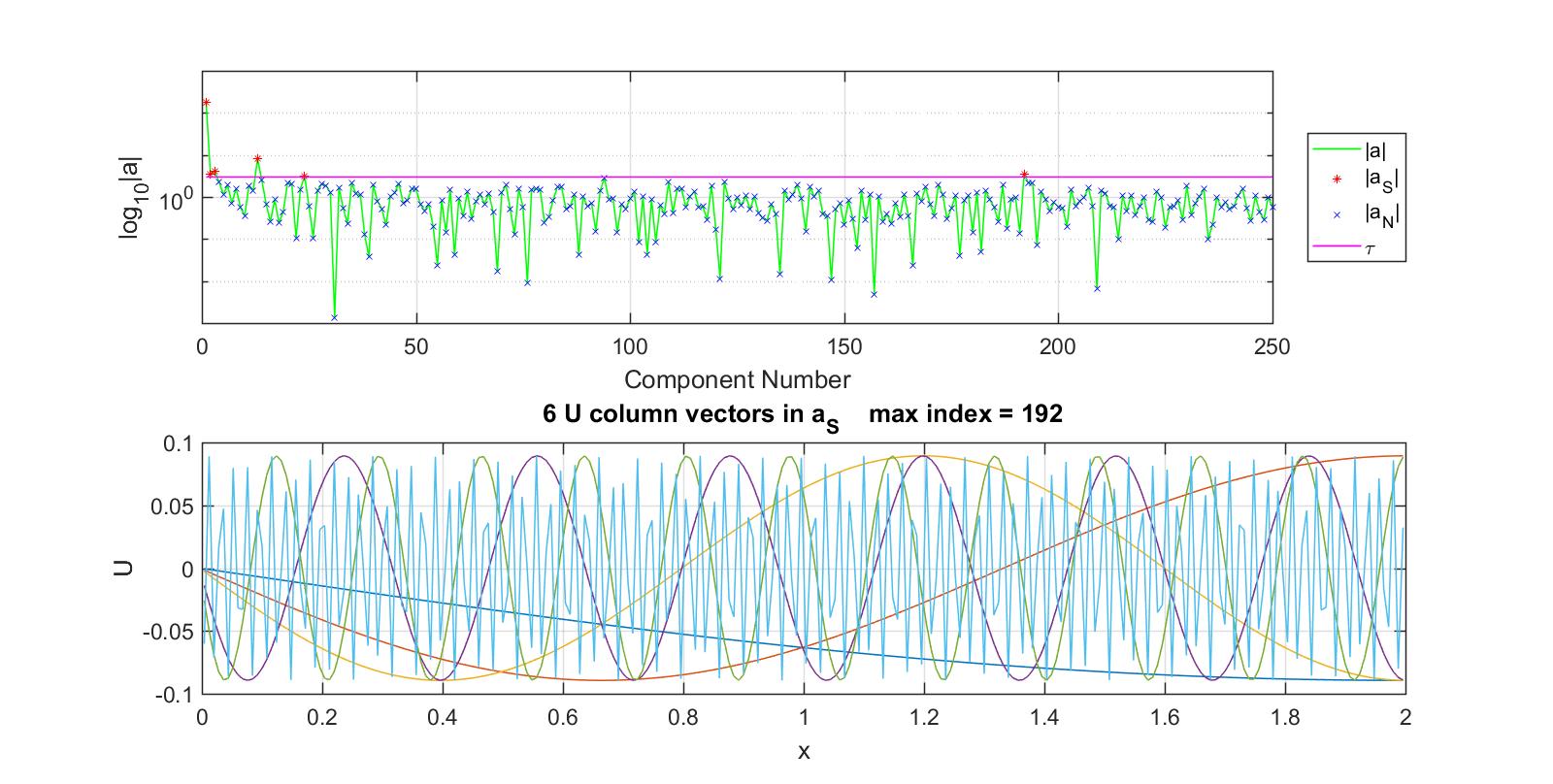}
  \caption{Upper:\ Separation of signal from noise in rotated data vector $\mathbf{a}$. Lower:\ Plot of colunms of $\mathbf{U}$ that correspond to signal components in $\mathbf{a}$.}
  \label{fig:aprojMidPt}
\end{figure}

\begin{figure}[tbp] 
  \centering
  \includegraphics[width=5in,height=2.52in,keepaspectratio]{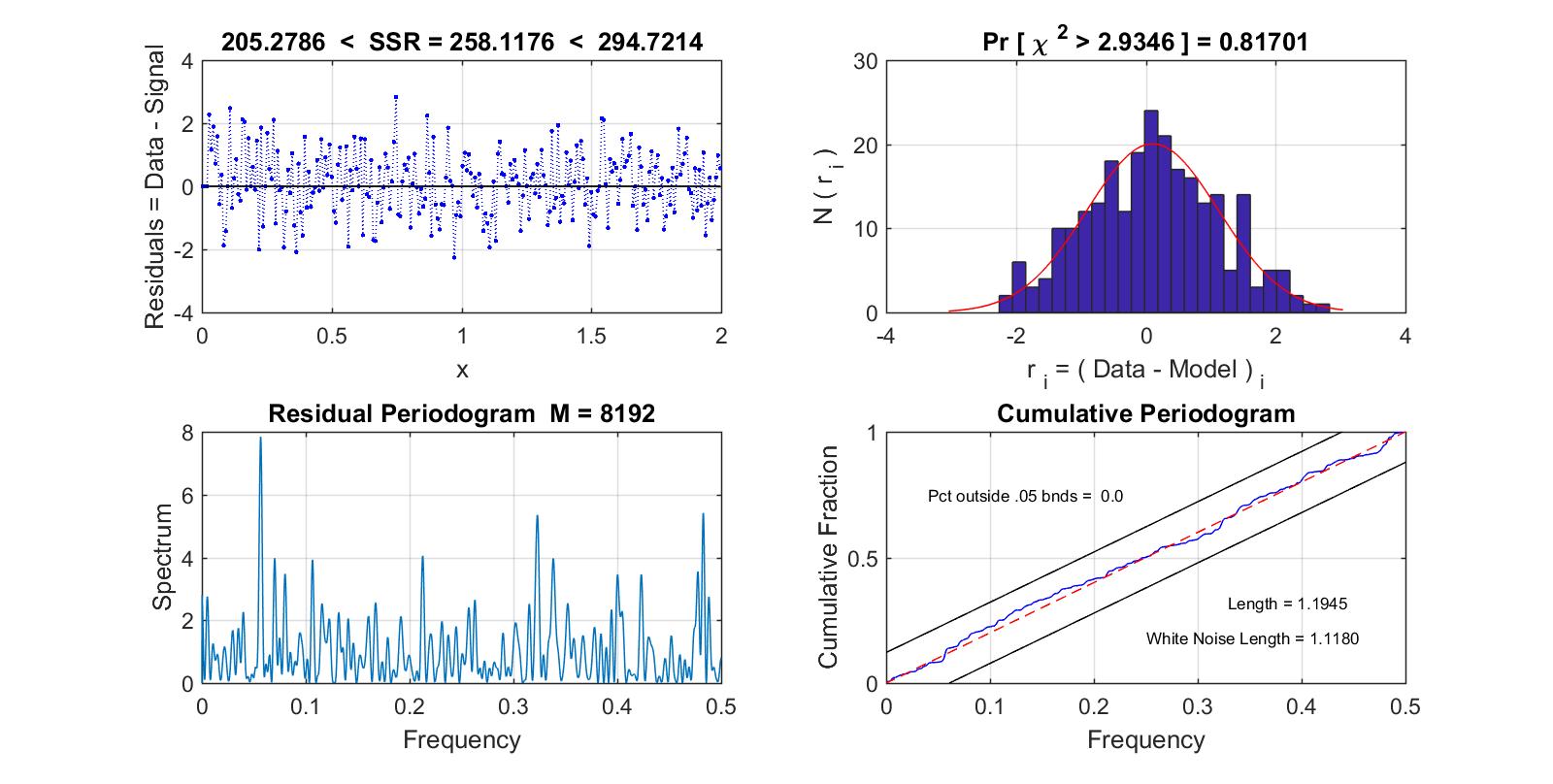}
  \caption{Diagnostics of residual vector, with $a_{192}$ included in signal vector $\mathbf{g}_S$. In top left figure, SSR stands for sum of squared residuals; upper and lower bounds on SSR correspond to $m \pm 2\sqrt{2m}$ bounds in \req{roughbounds}. Information above top right set of plots is result of chi-square goodness of fit test for null hypothesis that data in 
$\hat{\mathbf{r}}$ come from a normal distribution with mean and variance estimated from $\hat{\mathbf{r}}$, at $5\%$ significance level; if $\mbox{Pr}\le 0.05$, then reject null hypothesis. $M$ above periodogram plot in lower left refers to number of terms in DFT \req{DFT}. In lower right, cumulative periodogram \req{cumper} is plotted (solid line), together with cumulative periodogram of pure white noise (dashed line), and two parallel lines that define $95\%$ confidence band for white noise.}
  \label{fig:midpt_stats_w192}
\end{figure}

Diagnostics $1$-$3$ are summarized in Fig.~\ref{fig:midpt_stats_w192}.
The three diagnostics indicate that we have found an acceptable truncation of $\mathbf{a}$. However, if we just go ahead and estimate the source with high-frequency component ${a}_{192}$ included in the signal $\mathbf{g}_S$, we get the results in 
Fig.~\ref{fig:Craig_gtilde_fhat_w_a192}, which indicate that ${a}_{192}$ should have been included in the noise.
In fact, in a series of measurements of a standard, normally distributed quantity, a measurement that exceeds three standard deviations will occur roughly once in every $370$ terms, which means that the occurrence of one such measurement amongst the higher-frequency components of 
$\mathbf{a}$ is not an unlikely event.
If we include ${a}_{192}$ in $\mathbf{a}_N$, we get the much more satisfactory results in Fig.~\ref{fig:midpt_gs_fhat_w24}, with the sum of squared residuals now equal to $270.9$, which still satisfies \req{roughbounds}. 
This leads us to add a fourth diagnostic for separating signal from noise in the data.

\vspace{1em}
\noindent \textbf{Diagnostic 4}. {\em Even if a higher-frequency component of the projected signal $\mathbf{a}$ satisfies $|{a}_j| > \tau$, include that component in the noise $\mathbf{a}_N$}. 
\vspace{1em}

Examining the estimate of the source function $f$ which corresponds to the regularized data function $g$ approximated by the linear combination of columns of $\mathbf{U}$, we note that there is a dominant frequency of oscillation, corresponding to $a_{13}$, together with an oscillation corresponding to $a_{24}$, at almost twice the frequency. Based on the smoothing property \req{filter} of the convolution kernel \req{kernel}, we decide also to include ${a}_{24}$ in the noise. Alternatively, we could get the same result by adjusting $\tau$ upward to $3.2$. Including $a_{24}$ in $\mathbf{a}_N$, we get that the sum of squared residuals is equal to $280.6$, which still satisfies \req{roughbounds}. The final results are plotted in Fig.~\ref{fig:midpt_gs_fhat_w13} and \ref{fig:midpt_stats_w13}.

\begin{figure}[tbp] 
  \centering
  \includegraphics[width=5in,height=2.89in,keepaspectratio]{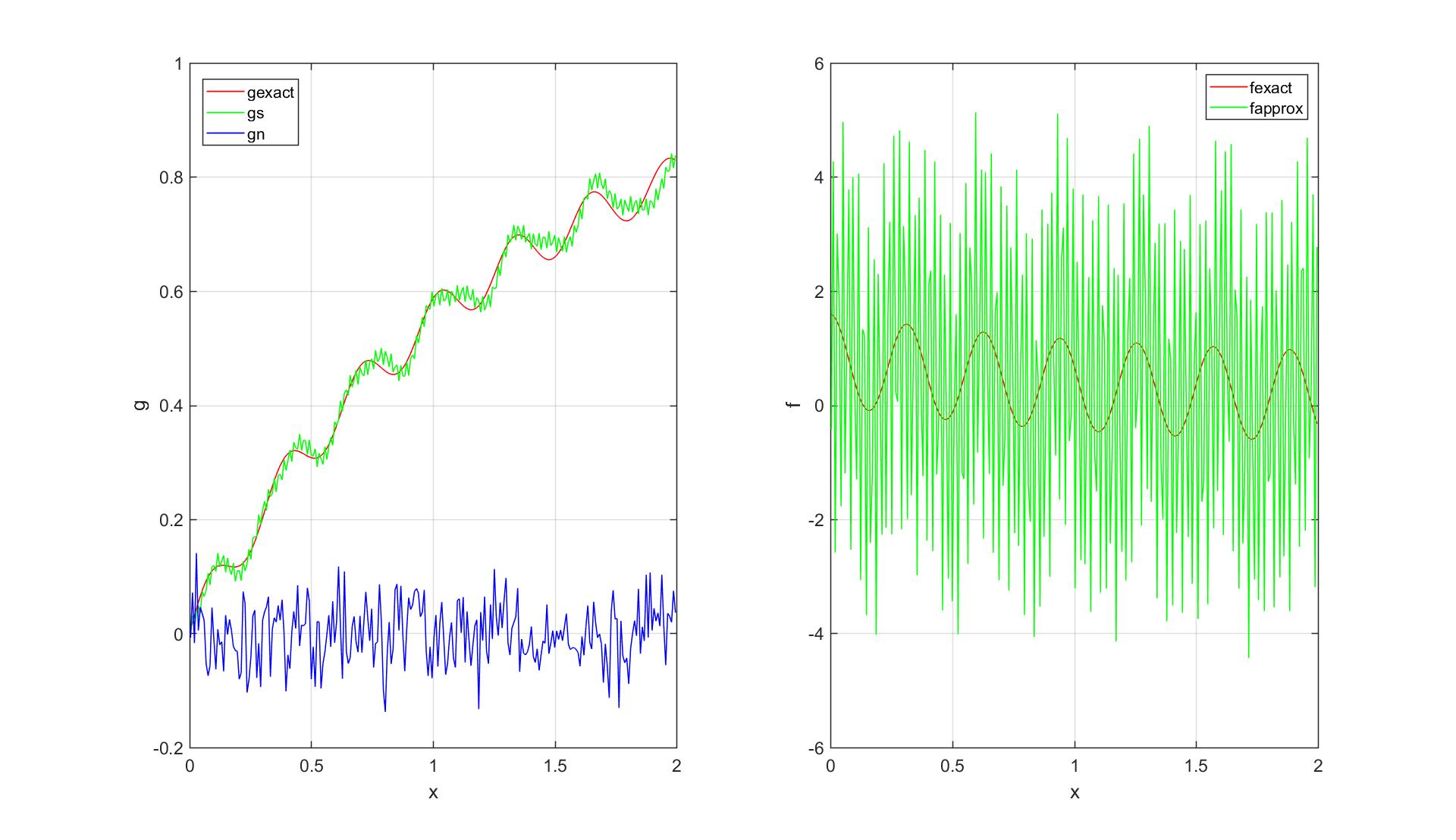}
  \caption{Plots of truncated data $g$ and corresponding source function $f$ with $a_{192}$ included in signal vector $\mathbf{g}_S$.}
   \label{fig:Craig_gtilde_fhat_w_a192}
\end{figure}

\begin{figure}[tbp] 
  \centering
  \includegraphics[width=5in,height=2.52in,keepaspectratio]{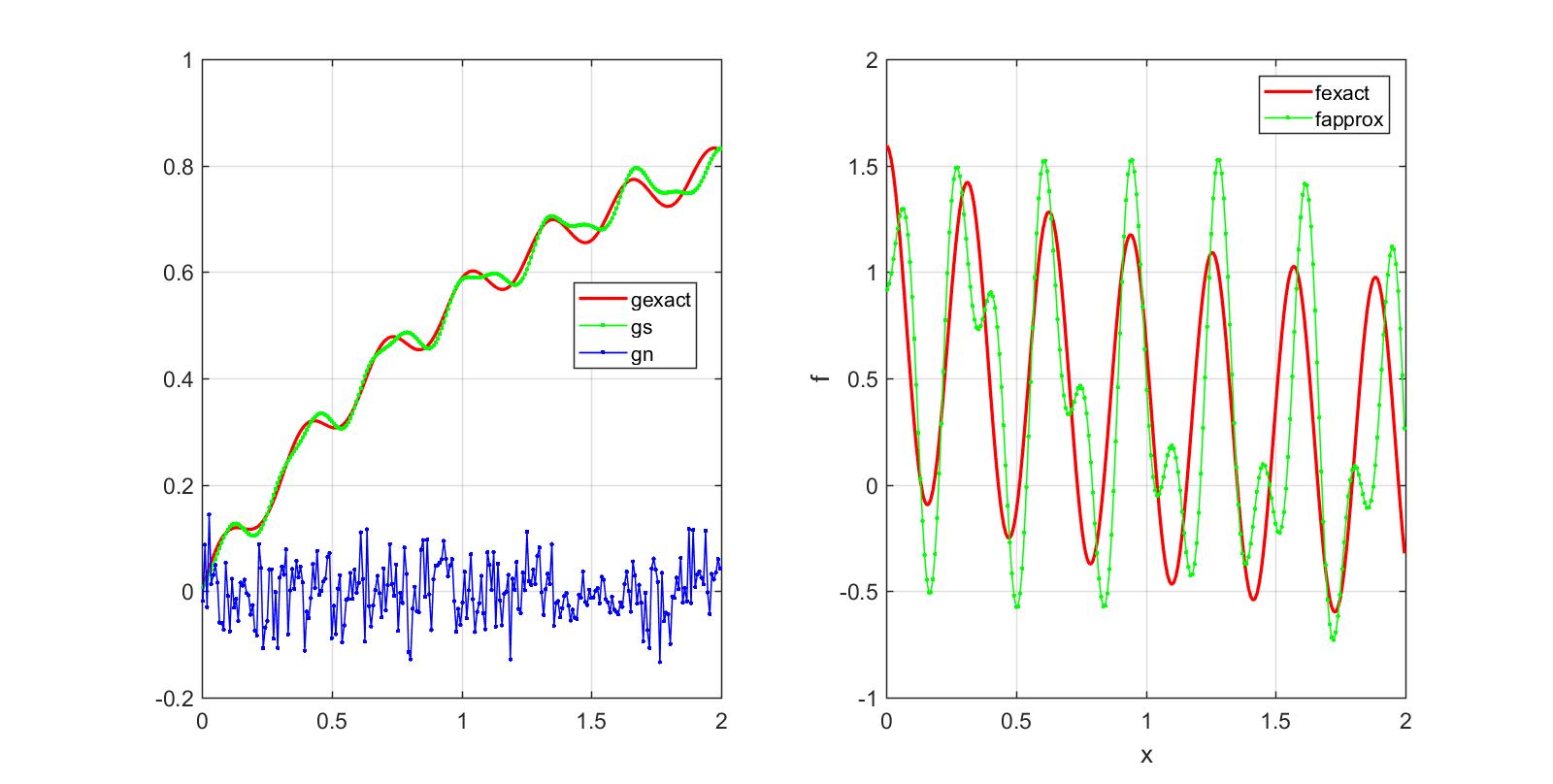}
  \caption{Plots of truncated data $g$ and corresponding source function $f$ with $a_{192}$ included in noise vector $\mathbf{g}_N$.}
  \label{fig:midpt_gs_fhat_w24}
\end{figure}

In the next section, we will discuss a related projection method for separating signal from noise in a discrete set of noisy data, that is based on closed-form SVD's of the Hilbert space operators of integration \req{integration} and fractional integration \req{Abel}.

\begin{figure}[tbp] 
  \centering
  \includegraphics[width=5in,height=2.52in,keepaspectratio]{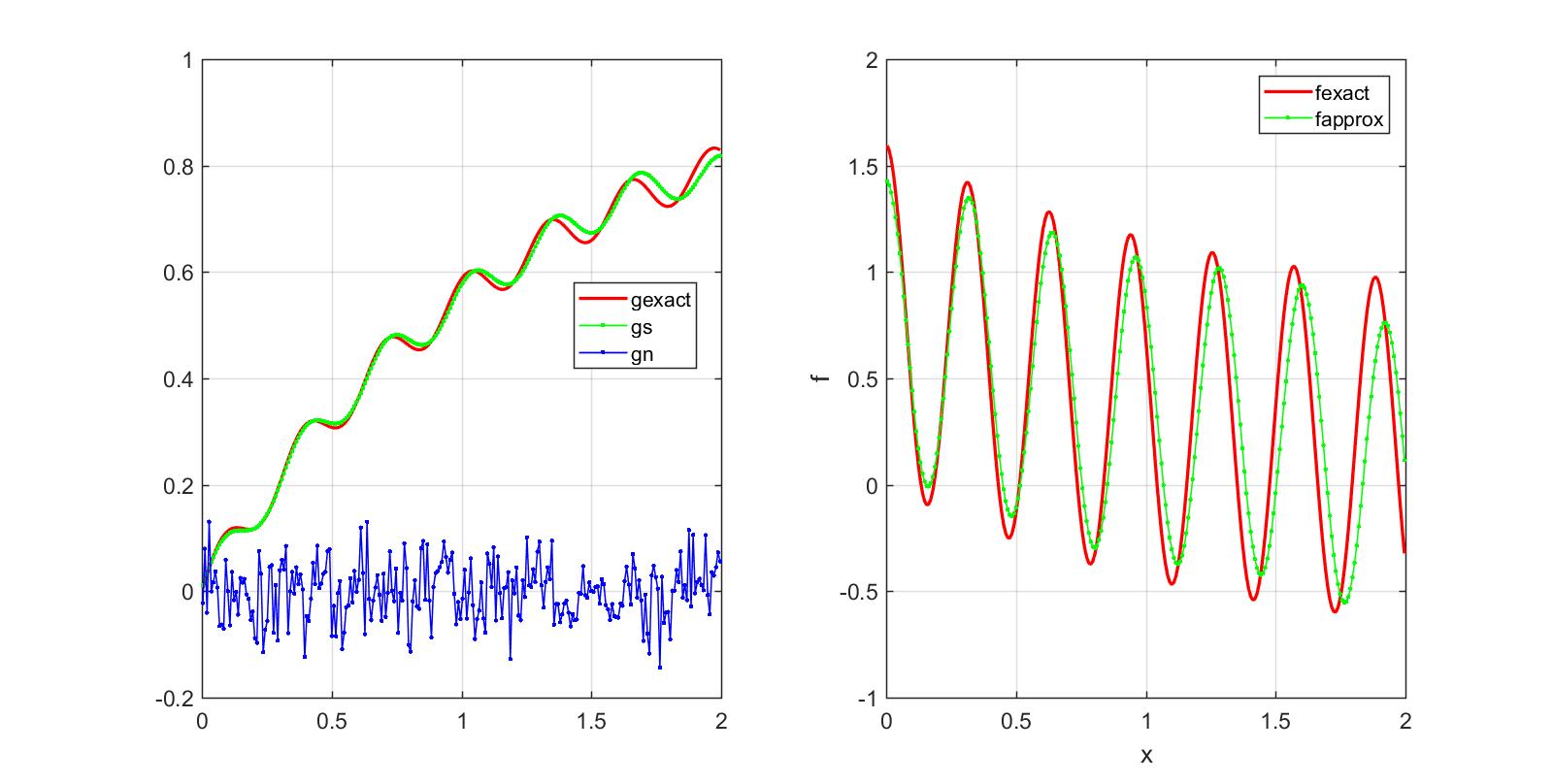}
  \caption{Plots of truncated data $g$ and corresponding source function $f$ with $a_{24}$ also included in noise vector $\mathbf{g}_N$.}
  \label{fig:midpt_gs_fhat_w13}
\end{figure}

\begin{figure}[tbp] 
  \centering
  \includegraphics[width=5in,height=2.9in,keepaspectratio]{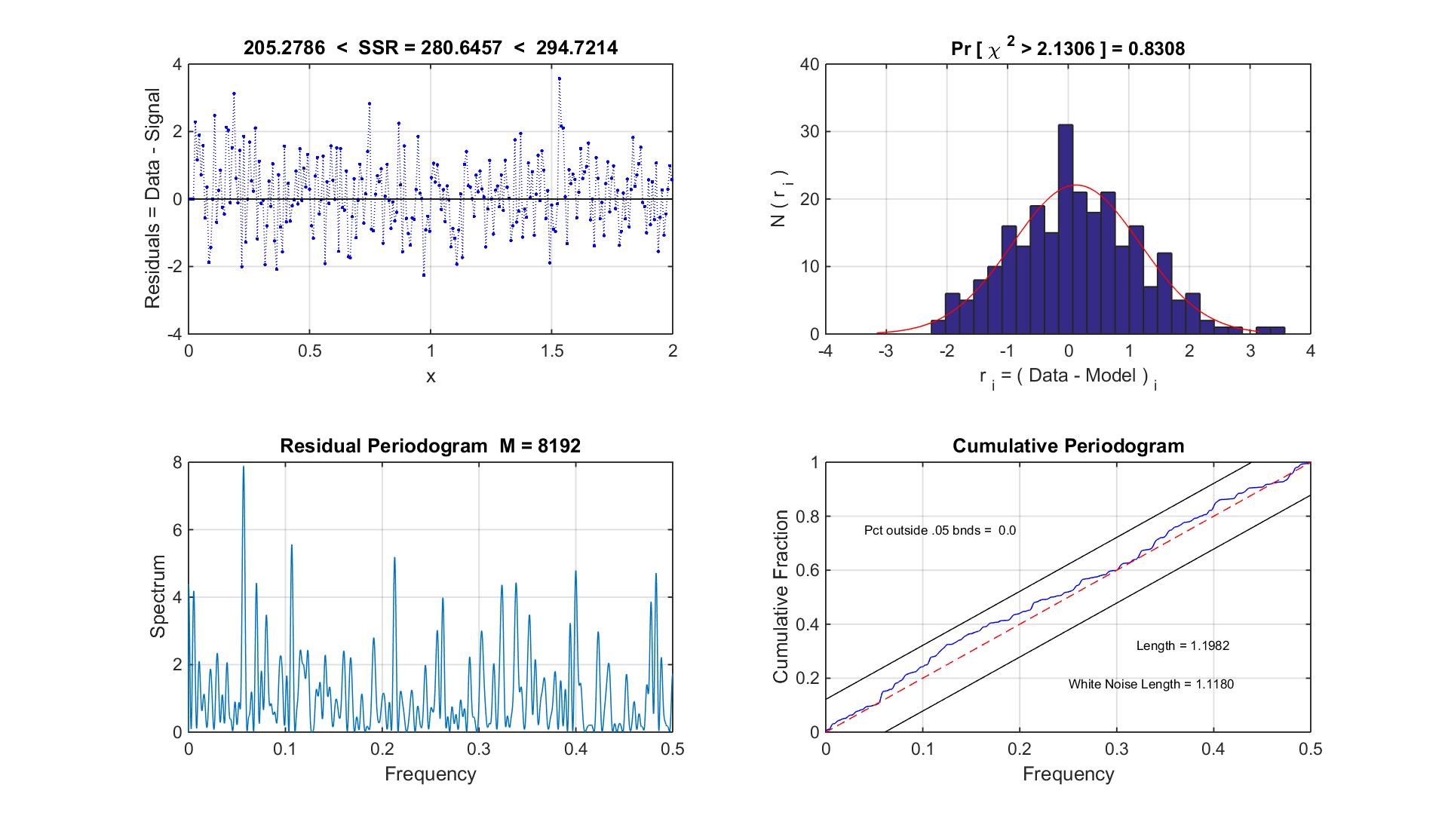}
  \caption{Diagnostics of residual with $a_{24}$ and $a_{192}$ included in noise.}
  \label{fig:midpt_stats_w13} 
\end{figure}

\section{New Projection Method}\label{sec:projection}
In Rust's method of regularization, as outlined in the preceding section, the orthonormal column vectors $\{\mathbf{u}_j,\ j=1,\ldots,m\}$ of the matrix $\mathbf{U}$, which are obtained as part of the discrete SVD \req{finSVD}, depend upon the discretization of the infinite-dimensional problem $Af=g$ that produces the finite-dimensional full-rank regression matrix $\mathbf{M}$. These vectors are ordered by the associated singular values of $\mathbf{M}$, $(\sigma_1,$ $\ldots,$ $\sigma_m)$, which are listed in order of decreasing size. In this section, for the operators of fractional integration \req{Abel} and integration \req{integration}, we will present a different projection method, that will provide a finite-dimensional approximation of $Af=g$, in which the oscillatory behavior of the basis vectors of both the domain and the range of the associated finite-dimensional operator can be characterized in detail. The method is based on the infinite-dimensional version of the SVD for compact linear operators on Hilbert space.

\subsection{Singular Value Decomposition}\label{subsec:infSVD1}
Let $X$ and $Y$ be infinite-dimensional real Hilbert spaces, with respective inner products $(\cdot\,,\cdot)_X$ and $(\cdot\,,\cdot)_Y$, and respective norms $\|\,\cdot\,\|_X$ and $\|\,\cdot\,\|_Y$ defined by the inner products. Let
$A: X\rightarrow Y\labeq{AXY}$
 be an injective compact linear operator, and let $A^\ast : Y\rightarrow X$ be its adjoint. It can be shown
 \citep[Thm.\ 15.16]{kress} that there exists a {\em singular system} 
 $\{(\sigma_j,\ v_j,\ u_j),\ j=1,2,\ldots\}$ of the operator $A$, consisting of a nonincreasing sequence of positive numbers $\{\sigma_j,\ j=1,2,\ldots\}$ with limit $0$ as $j\rightarrow\infty$, an orthonormal basis $\{v_j,\ j=1,2,\ldots\}\subset X$, and a corresponding orthonormal system $\{u_j,\ j=1,2,\ldots\}\subset Y$, such that
\be 
A v_j = \sigma_j u_j,\hspace{3ex}A^\ast u_j = \sigma_j v_j,\hspace{3ex}A^\ast A v_j = \sigma^2_j v_j,\hspace{3ex} j=1,2,\ldots\,.\labeq{Avj}\ee
In addition, for each $f\in X$, we have the {\em Singular Value Decomposition}
\be 
f=\sum_{j=1}^\infty \left(f,v_j\right)_X v_j,\hspace{5ex}\mbox{and}\hspace{5ex}Af =\sum_{j=1}^\infty \left(f,v_j\right)_X \sigma_j u_j.\labeq{Af}\ee
Following \citet{cuer}, the SVD can be described in ``matrix'' notation as follows. The action of the operator $A$ is a composition of three mappings,
\be 
A = U\Sigma V^T,\labeq{svd}\ee
such that 
\be 
V^T : X\rightarrow \ell^{\,2}, \hspace{3ex}\Sigma :\ell^{\,2}\rightarrow \ell^{\,2},\hspace{3ex}U:\ell^{\,2}\rightarrow Y,\labeq{USVT}\ee
where $\ell^{\,2}$ is the space of real square-summable sequences, $\Sigma=\mbox{diag}\left(\sigma_j\right)$, and 
\bea
 V^T &:&f\ \mapsto\  \{\left(f,v_j\right)_X,\ j=1,2,\ldots\,\},\nonumber\\[4pt]
 \Sigma^{\ } &:&\{\left(f,v_j\right)_X,\ j=1,2,\ldots\,\}\ \mapsto\  \{\sigma_j\left(f,v_j\right)_X,\ j=1,2,\ldots\,\},\nonumber\\[4pt]
 U^{\ } &:&\{\sigma_j\left(f,v_j\right)_X,\ j=1,2,\ldots\,\}\ \mapsto\   \sum\nolimits_{j=1}^\infty \sigma_j\left(f,v_j\right)_X u_j.\nonumber\eea
The following theorem is due to Picard \citep[Thm.\ 15.18]{kress}.
\begin{theorem*} 
The equation of the first kind 
\be Af=g\labeq{Afg}\ee 
is solvable if and only if $g\in Y$, and
\be 
\sum_{j=1}^\infty { \left( \frac{|\left(g,u_j\right)_Y|}{\sigma_j}\right) }^2 < \infty\labeq{Picard}.\ee
For each $g\in Y$ satisfying \req{Picard}, the unique solution $f\in X$ of $Af=g$ is given by
\be 
f =\sum_{j=1}^\infty \frac{1}{\sigma_j}\left(g,u_j\right)_Y v_j.\labeq{fsol}\ee
\end{theorem*}

Note that the {\em Picard condition} \req{Picard} excludes many functions $g\in Y$ from $A(X)$, the range of $A$.
Suppose that $f$ is the solution of \req{Afg} for some $g$ satisfying \req{Picard}. If we add a small perturbation to $g$, $g^\delta = g+\delta u_j$, we get that $f^\delta = f + \delta\sigma_j^{-1}\,v_j$. It follows from \req{fsol} that 
${\|f^\delta -f\|_X}/{\|g^\delta -g\|_Y}={1}/{\sigma_j}$
can be made arbitrarily large, because the singular values decay to zero as $j\rightarrow\infty$. This demonstrates that equation \req{Afg} is ill-posed, and it is more ill-posed the faster the singular values of $A$ approach zero.

\subsection{Projection Method for Integration}\label{subsec:integration}

It turns out that, for $X = Y = L^2[0,1]$, the singular system $\{(\sigma_j,$ $v_j,$ $u_j),$ $\ j=1,2,\ldots\}$ for the operator of integration $A := I_0^1 $ \req{integration} has been determined \citep[Example 15.19]{kress}. It is given by
\begin{equation}
\sigma_j = \frac{1}{\pi c_j },\hspace{3ex}v_j(x)= \sqrt2 \cos(c_j\pi x),\hspace{3ex}u_j(x)= \sqrt2 \sin(c_j\pi x),\labeq{svdint}\ee
where
\be  
c_j = j-\nicefrac{1}{2},\hspace{5ex} j=1,2,\ldots\,,\labeq{cj}\ee
so that $\sigma_j\sim 1/j$ as $j\rightarrow\infty.$
In this case, the SVD of $A$ suggests a natural method for obtaining an approximate solution of $Af=g$. 

Define the finite-dimensional subspaces $X_m\subset X$ and $Y_m\subset Y$, as follows,\be X_m = \mbox{span}\{v_j:\ j=1,\ldots,m\},\hspace{5ex} 
    Y_m = \mbox{span}\{u_j:\ j=1,\ldots,m\}.\labeq{XmYm}\ee
As in Sec.~\ref{sec:rustreg}, suppose that the smooth function $g$ is defined by the noisy data set $\{g_1,\ldots,g_m\}$, corresponding to $0 < x_1 < x_2 < \cdots <x_m =1$, and assume that the subspace $Y_m$ is unisolvent with respect to these points, i.e., each function in $Y_m$ that vanishes at $x_1,$ $\ldots,$ $x_m$ vanishes identically. This is certainly true for the $m$ basis functions $u_1,\ldots,u_m$ in $ Y_m$ in \req{svdint}, so that it is also true for every function in $Y_m$. 
Then there exists a unique function 
\be  G(x)\ =\ \sum_{j=1}^m \xi_j u_j(x),\hspace{2ex}k=1,\ldots,m,\labeq{G}\ee which is infinitely differentiable on $[0,1]$, satisfies $G(0)=0$, and  interpolates the given data,
\be g_k =  G(x_k)\ =\ \sum_{j=1}^m \xi_j u_j(x_k),\hspace{2ex}k=1,\ldots,m.\labeq{gminterp}\ee
Since $G(x)$ is a finite sum, it satisfies \req{Picard}, so that $G\in A(X)$.

The unique set of $m$ real coefficients $\boldsymbol{\xi}$ 
can be determined by solving the matrix equation
\be \mathbf{g} = \mathbf{P}\boldsymbol{\xi},\labeq{gBxi}\ee
where
\be
\mathbf{P} = 
\begin{pmatrix}
    u_1(x_1) & u_2(x_1) & \dots  & u_m(x_1) \\
    u_1(x_2) & u_2(x_2) & \dots  & u_m(x_2) \\
    \vdots & \vdots & \ddots & \vdots \\
    u_1(x_m) & u_2(x_m) & \dots  & u_m(x_m) \\
\end{pmatrix}
,\labeq{Upsmatrix}\ee
so that
\be \boldsymbol{\xi}=\mathbf{P}^{-1}\mathbf{g}.\labeq{xisol}\ee
It follows from \req{fsol} that an estimate $\hat{f} \in X_m$ of the source function $f\in X$ is given by the finite-dimensional linear combination,
\be 
\hat{f}(x) =  \sum_{j=1}^m \frac{\xi_j}{\sigma_j}\, v_j(x), \labeq{hatfgen}\ee
which, for the operator of integration \req{integration}, can be written as follows,
\be 
\hat{f}(x)  = \sqrt2\,\pi \sum_{j=1}^m \xi_j  c_j \cos(c_j\pi x),\labeq{hatfint}\ee
where $c_j = j-\nicefrac12.$

We are once again faced with the same kind of instability problem that was discussed in Section \ref{sec:rustreg}. Because of the ill-posedness of the underlying integral equation, without prior regularization of the data $\mathbf{g}=(g_1,\ldots,g_m)^T$, the estimate \req{hatfint} for the source function $f(x)$ would be overwhelmed by amplified high-frequency noise in the data. This amplification results from division of the components in $\boldsymbol{\xi}$ by successively smaller singular values in each case. We will discuss how to resolve these difficulties in Sec.~\ref{sec:newreg}. First, however, we want to show how we can use a similar procedure based on the SVD to obtain an estimate for the source function in the case of fractional integration.

\subsection{Projection Method for Fractional Integration}\label{subsec:fracint}
It is apparently much less well-known that, for $\mu\in(0,1)$,  $\Omega=[-1,1]$, $X=L^2(\Omega,w_1)$, and
$Y=L^2(\Omega,w_2)$, \citet{GorenfloTuan1995} have determined a family of singular systems $\{(\sigma_j,\ v_j,\ u_j),\ j=1,2,\ldots\}$ for the Abel transform $I_{-1}^{\mu}$ \req{Abel}. In particular, they have shown that, with $w_1(x)=1$ and $w_2(x)= (1-x^2)^{-\mu}$, the singular system is given by
\be  
\sigma_j = \sqrt{ \frac{\Gamma(j-\nicefrac{1}{2}-\mu)}{\Gamma(j-\nicefrac{1}{2}+\mu)} }\sim {1}/{j^{\mu}},\hspace{5ex}j\rightarrow\infty,\labeq{sigmamu}\ee
\be  
v_j(x) = \sqrt{2j-1}\,C_{j-1}^{\nicefrac{1}{2}} (x),\hspace{5ex}u_j(x) = c_j\,(1+x)^\mu\, P_{j-1}^{(-\mu,\mu)}(x),\labeq{uvmu}\ee
where the $C_{j-1}^{\nicefrac{1}{2}}$ are Gegenbauer polynomials of degree $j-1$, $P_{j-1}^{(-\mu,\mu)}$ are Jacobi  polynomials of degree $j-1,$ and 
\be  
c_j = \sqrt{ \frac{(j-\nicefrac{1}{2})\,\left[\Gamma(j)\right]^2 }{\Gamma(j-\mu)\Gamma(j+\mu)} }.\labeq{cjmu}\ee
We note that 
\be C_j^{\nicefrac{1}{2}}(x) = P_j(x),\hspace{5ex} j=0,1,\ldots,\labeq{gegenauer+legendre}\ee 
where $P_j$ are Legendre polynomials \citep[18.7.9]{NIST:DLMF}, and that $u_j(x)$ is not differentiable at $x=-1$, $j=1,2,\ldots$.

As before,
we determine the unique function $G(x)\in Y_m\subset L^2(\Omega,w_2)$ that interpolates the data $\mathbf{g}$, by solving for the unique set of coefficients $\xi_1,\ldots,\xi_m$, such that
\be g_k  = \sum_{j=1}^m \xi_j (1+x_k)^\mu P_{j-1}^{(-\mu,\mu)}(x_k),\hspace{2ex}k=1,\ldots,m.\labeq{gmJacobi}\ee
This function satisfies $G(-1)=0$, and is smooth on $(-1,1]$.
To do this, we form the analogue of the matrix $\mathbf{P}$ in \req{Upsmatrix}, 
\be \mathbf{P}= 
\begin{pmatrix}
    u_1(x_1) & u_2(x_1) & \dots  & u_m(x_1) \\
    u_1(x_2) & u_2(x_2) & \dots  & u_m(x_2) \\
    \vdots & \vdots & \ddots & \vdots \\
    u_1(x_m) & u_2(x_m) & \dots  & u_m(x_m) \\
\end{pmatrix},\labeq{JacobiC}\ee
where $u_j(x) = (1+x)^\mu P_{j-1}^{(-\mu,\mu)}(x),\ j=1,\ldots,m$, and get 
\be \boldsymbol{\xi}=\mathbf{P}^{-1}\mathbf{g}.\labeq{xisol2}\ee
For fractional integration, this provides an estimate of the source function $f$, given by the finite linear combination of Legendre polynomials,
\be 
\hat{f}(x) = \sum_{j=1}^m \frac{\xi_j}{\sigma_j}\, \sqrt{2j-1}\,P_{j-1} (x).\labeq{hatffracint}\ee
We must still address the problem of regularizing the noisy data $\mathbf{g}$ prior to obtaining an estimate of ${f}$. This will be discussed in the next section.

\section{Regularization by Truncated Projection} \label{sec:newreg}

Our initial goal in this section is to develop a method for separating signal from noise in the data, $\mathbf{g}=\mathbf{g}_S+\mathbf{g}_N$, which is similar to Rust's method in Sec.~\ref{sec:rustreg}, and which also takes advantage of the fact that the SVD's for integration and fractional integration are known explicitly. We will then show how an extension of this method can be used to regularize discrete sets of data measurements which have been contaminated by random, zero-mean noise, for other applications. 

\subsection{Regularization Using Closed-Form SVD}\label{subsec:infSVD}

As before, we first scale each data point by its estimated standard deviation, and get the scaled data vector $\mathbf{b}=\mathbf{S}^{-1}\mathbf{g}$. We want to regularize $\mathbf{b}$, i.e., we to compute a decomposition $\mathbf{b}=\mathbf{b}_S+\mathbf{b}_N$, such that $\mathbf{b}_N $ captures a sufficient amount of the variance in the data, i.e., such that $\|\mathbf{b}_N\|_2^2$ lies within the bounds \req{roughbounds}.
Rather than do this by computing the discrete singular value decomposition of the matrix $\mathbf{P}$ in Sec.~\ref{subsec:integration} or \ref{subsec:fracint} as we did with the matrix $\mathbf{M}$ in Sec.~\ref{sec:rustreg}, our approach here is to compute instead its $QR$ decomposition, 
\be\mathbf{P}=\mathbf{Q}\mathbf{R}.\labeq{PQR}\ee
Here, $\mathbf{Q}$ is an $m\times m$ orthogonal matrix whose columns span $\mathbb{R}^m$, and $\mathbf{R}$ is an invertible $m\times m$ upper triangular matrix. 
We then project $\mathbf{b}$ onto the successive orthonormal columns of $\mathbf{Q}$, and get
\be \mathbf{a} =\mathbf{Q}^T\mathbf{b} .\labeq{wQTb}\ee

Because $\mathbf{Q}$ is an orthogonal matrix, $\|\mathbf{a}\|_2 = \|\mathbf{b}\|_2$.
Since the data in vector $\mathbf{b}$ have been scaled to units of one standard deviation, the same is true of the data in the vector $\mathbf{a}$. We proceed as in Sec.~\ref{sec:rustreg}, 
and separate $\mathbf{a}$ into signal and noise, respectively, according to whether or not a component exceeds a specified truncation level, $|a_k|>\tau$, starting with our rule of thumb, $\tau = 3$. As we will demonstrate in the next section, 
\be k_{\max} = \max\{k:|a_k|>\tau,\ \ k=1,\ldots,m\}\labeq{maxk}\ee 
is typically significantly smaller than the number of data points $m$.

We thus obtain the decomposition 
\be\mathbf{a}=\mathbf{a}_S + \mathbf{a}_N,\labeq{atrunc}\ee
 so that we also have the decomposition of the scaled data,
\be \mathbf{b}=\mathbf{b}_S+\mathbf{b}_N = \mathbf{Q}\mathbf{a}_S+\mathbf{Q}\mathbf{a}_N,\labeq{btrunc}\ee
and of the original data,
\be \mathbf{g}=\mathbf{g}_S+\mathbf{g}_N = \mathbf{S}\mathbf{b}_S+\mathbf{S}\mathbf{b}_N.\labeq{gtrunc}\ee
We then proceed to construct an estimate \req{hatfgen} of the source function, that is based only on $\mathbf{g}_S$. 
Let 
\be \boldsymbol{\xi} = \mathbf{P}^{-1}\mathbf{g}_S = \mathbf{R}^{-1}\mathbf{Q}^{T}\mathbf{g}_S .\labeq{Pinvsig}\ee
A closed-form estimate $\hat{f}(x)$ of the source function $f(x)$ is then given by either \req{hatfgen} or \req{hatffracint}. 

\subsection{Frequency Content of the Data}\label{subsec:frequency}

There is an important consequence of using the $QR$ factorization \req{PQR} to separate signal from noise in the scaled, rotated data $\mathbf{a}$ \req{wQTb}, and to determine the coefficients $\boldsymbol{\xi}$ when a closed-form singular system $\{(\sigma_j,$ $v_j,$ $u_j),$ $\ j=1,2,\ldots\}$ is available for the operator $A$, and the functions  
$v_j$ and $u_j$ oscillate like  orthogonal polynomials of degree $j-1$ on $(a,b)$. In this case, each of $v_j$ and $u_j$ has exactly $j-1$ zeros in $a<x<b$, so that the number of oscillations of $v_j$ and $u_j$ in $(a,b)$ increases as $j$ increases.
As long as we use a $QR$ decomposition algorithm that does not permute the columns of the original matrix $\mathbf{P}$, it follows that, for $k>1$, each component $a_k$ encodes higher frequency content of the data $\mathbf{g}$ than the preceding component $a_{k-1}$, as $k$ increases. As we will demonstrate in the next section, just as in the example that was presented at the end of Sec.~\ref{sec:rustreg}, $\mathbf{a}_S$, the signal portion of $\mathbf{a}$, is determined by a relatively small subset of the components of $\mathbf{a}$, which consists of a few lower-index, lower-frequency components. Since $\mathbf{b}_S = \mathbf{Q}\mathbf{a}_S$, and $\mathbf{g}_S = \mathbf{S}\mathbf{b}_S$, it will follow that the smooth, closed-form function $G(x)\in Y_m$  that approximates the data function $g(x)$ is formed by a linear combination of only a few lower-index basis functions $ u_j(x) $.

First, however, we want to finish this section by introducing a more general method for regularizing noisy data, that is based on the approach that we have presented in Sec.~\ref{subsec:infSVD}, but which does not require the SVD of an operator. This method will provide a smooth, regularized polynomial approximation of a function that is defined by noisy data measurements, with estimated pointwise random errors, for other applications.

\subsection{General Method for Regularizing Noisy Data}\label{subsec:otherapps}
 
It is worth noting that, for the purpose of regularizing an unknown smooth function $g$, which is defined only by a finite set of noisy data, the usefulness of having a closed-form SVD for the operator that models the process that produced the data, is that it provides a set of basis functions, which span the range of the operator, and are easy to compute. Using such a basis, for the operators of integration \req{integration} and fractional integration \req{Abel}, we have shown in Sec.~\ref{subsec:infSVD} how to construct a smooth, closed-form function $G(x)$, which approximates the unknown data function $g$, in the norm of $L^2([a,b],w)$, for some weight function $w$. We can introduce a different set of smooth basis functions, which will ``decouple" our method from {\em any} SVD, as follows.

Suppose that we have been given a discrete set of $m$ data measurements $ \mathbf{g} = (g_1,$ $\ldots,g_m)^T$, corresponding to a subdivision (which does not need to be equally spaced) $a \le x_1<x_2<\cdots <x_m \le b$ of $[a,b]$, which represent a smooth function $g\in C^1[a,b]$. Without loss of generality, we assume that $[a,b]$ is the interval $[-1,1]$. Also suppose that the data have been contaminated by random errors 
$\boldsymbol{\epsilon}=(\epsilon_1,\ldots,\epsilon_m)^T$, for which it is known that the conditions \req{calE} and \req{S2} in Sec.~\ref{sec:rustreg} are satisfied. Finally, suppose that these data have been generated by some process that smoothes the high-frequency content of the source that produced our data.  We can obtain a regularized polynomial approximation of the data function $g$ by the following procedure. 
  
Choose a convenient ordered set of Jacobi polynomials $p_j$ of degree $j-1$, for $j=1,\ldots,m$, which are orthonormal in 
$Z={L^2([-1,1],w)}$, where $w$ is the appropriate weight function. Denote by $Z_m$ the subspace of $Z$ that is spanned by this set of polynomials. Form the appropriate matrix $\mathbf{P}$, as in (\ref{subsec:integration}) and (\ref{subsec:fracint}), and then compute its $QR$ factorization. Separate low-frequency signal from higher-frequency noise in the scaled data, by truncating the scaled data, exactly as we have done above, so that we have the decomposition  $\mathbf{g}=\mathbf{g}_S + \mathbf{g}_N$.
Then solve for the unique set of $k_{\max}$ coefficients $\boldsymbol{\xi}$ of the $(k_{\max}-1)$-degree polynomial ${G}(x)\in Z_m$ that interpolates the signal portion of the data,
$\boldsymbol{\xi} = \mathbf{P}^{-1}\mathbf{g}_S $. 

Without the aid of a singular value decomposition of either a discrete finite-dimensional or a compact infinite-dimensional operator, we have shown how to obtain a smooth, closed-form polynomial function $G(x)\in Z_m \subset L^2([-1,1],w)$, whose values are easy to compute, and which approximates $g(x)$ by interpolating the regularized data $\mathbf{g}_S$ on the subdivision $-1 \le x_1<x_2<\cdots< x_m \le 1$. As we will see in the next section, this polynomial is typically of much lower degree than $m-1$. Furthermore, this method for constructing $G(x)$ gives useful information about the frequency content of a given set of measurements. 
In the next section, we will present the results of some numerical simulations that demonstrate the usefulness of our method of data regularization, and we will also demonstrate the usefulness and simplicity of estimating the source function in \req{integration} and 
\req{Abel} by means of closed-form singular systems for the respective operators.

\section{Numerical Examples}\label{sec:numerical}

All of the calculations in this paper have been performed on a uniform grid, using the 64-bit Windows version of Matlab (R2017a), and Jacobi polynomials have been evaluated using Matlab software that has been developed by Burkardt \citep{fsu+codes}. The $QR$ decompositions have been performed using the Matlab version of Gander's modification of Rutishauser's algorithm $\mathbf{gramschmidt}$ \citep{gander}. 

Although our focus is on regularizing a function that has been defined by noisy data, we will first study the usefulness of the closed-form SVD expansions for obtaining the derivative of a differentiable function.

\subsection{Craig and Brown's Problem Revisited}\label{subsec:uncorrupted}

As our first example, we revisit Case (a) of the problem of Craig and Brown in Sec.~\ref{sec:rustreg}. 
We transform the independent variable $x$, so that the problem is specified on $[0,1]$, as follows,
\bea g(x) &=& 1-\exp(-2\alpha x)+\beta\sin(2\omega x),\labeq{gcraig2} \\[5pt]
 f(x) &=&  2\alpha\exp(-2\alpha x)+2\beta\omega\cos(2\omega x),\labeq{fcraig2}\vspace{2em}\eea
with the parameter values 
$\alpha = 0.8$, $\omega = 20$, and $\beta = 0.04$. 

As a check on differentiation using \req{svdint}, we evaluate $g$ in \req{gcraig2} on the equally-spaced mesh 
\be x_0=0,\,x_1=h,\,x_2=2h,\ldots,\,x_m=1,\labeq{xmesh}\ee
 where $h=1/m$, for $m=250$. Note that both $g(0)=0$ and each of the $u_j(0)=0$. Furthermore, each of the $v_j(1) = 0$, and neither $g$ nor $f$ is periodic on $[0,1]$. We interpolate $g$ at the points $x_1,\ldots,x_m$ using the first $N$ sine functions in \req{svdint}, and then we estimate $f$ by $\hat{f}$, using \req{hatfint} with the first $N$ cosine polynomials, for $N=15,\,25,\,50$. The corresponding differences between the exact solutions and the polynomial estimates are plotted in Fig.~\ref{fig:CraigBrownUnpertBeta04error}. There is nonuniform convergence of the approximation to $g$, and a Gibbs phenomenon in the approximation to $f$, at $x=1$. By $N=50$, the approximations are close, except at the right boundary of $f$; see Fig.~\ref{fig:CraigTrigUnpertbeta04N50}.

\begin{figure}[tbp] 
  \centering
  \includegraphics[width=5in,height=2.35in,keepaspectratio]{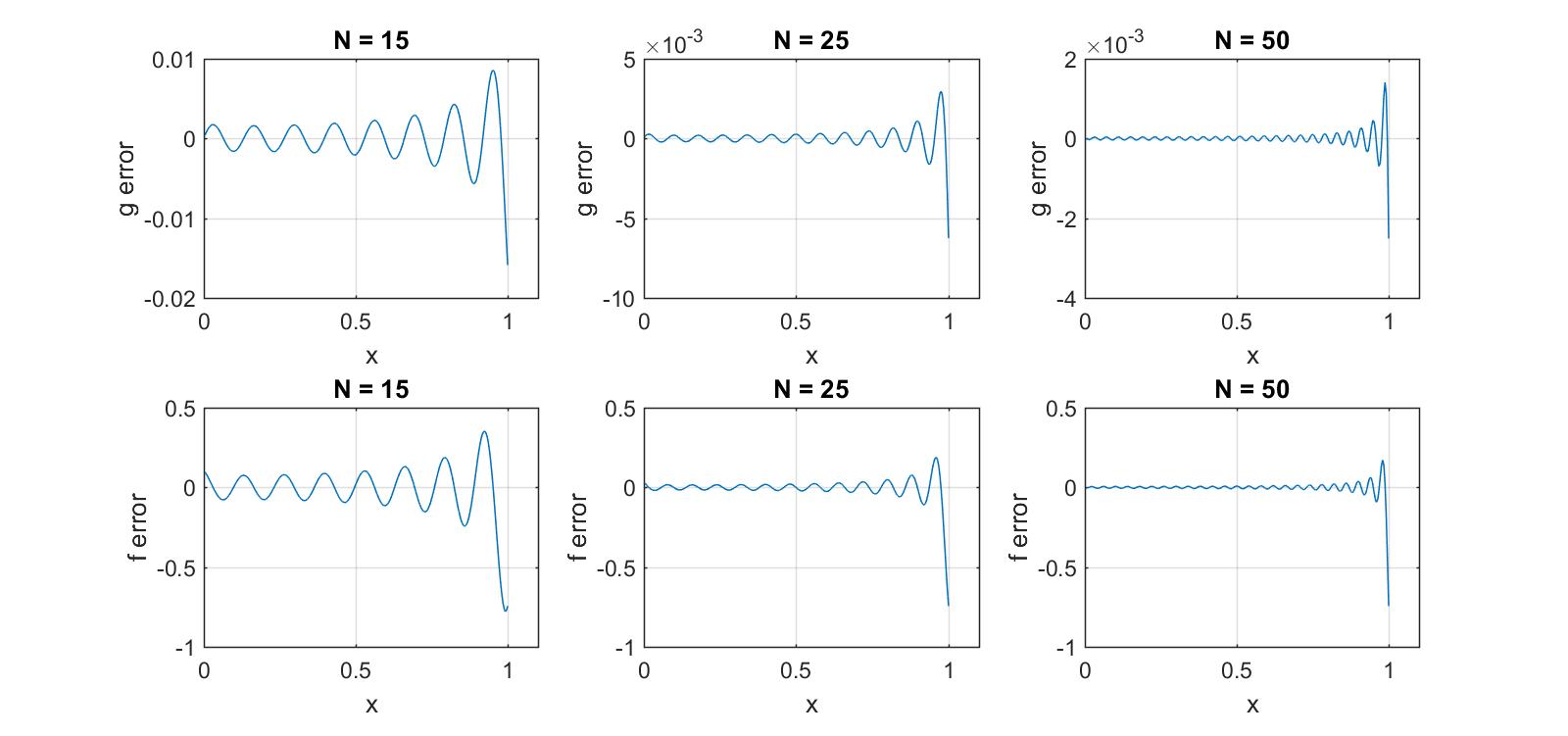}
  \caption{Differences between exact $g$ \req{gcraig2} (upper), and exact $f$ \req{fcraig2} (lower) and respective approximations of these functions using sums from $1,\ldots,N$ of trigonometric SVD expansions in Sec.~\ref{subsec:integration}.}
  \label{fig:CraigBrownUnpertBeta04error}
\end{figure}

\begin{figure}[tbp] 
  \centering
  \includegraphics[width=5in,height=2.2in,keepaspectratio]{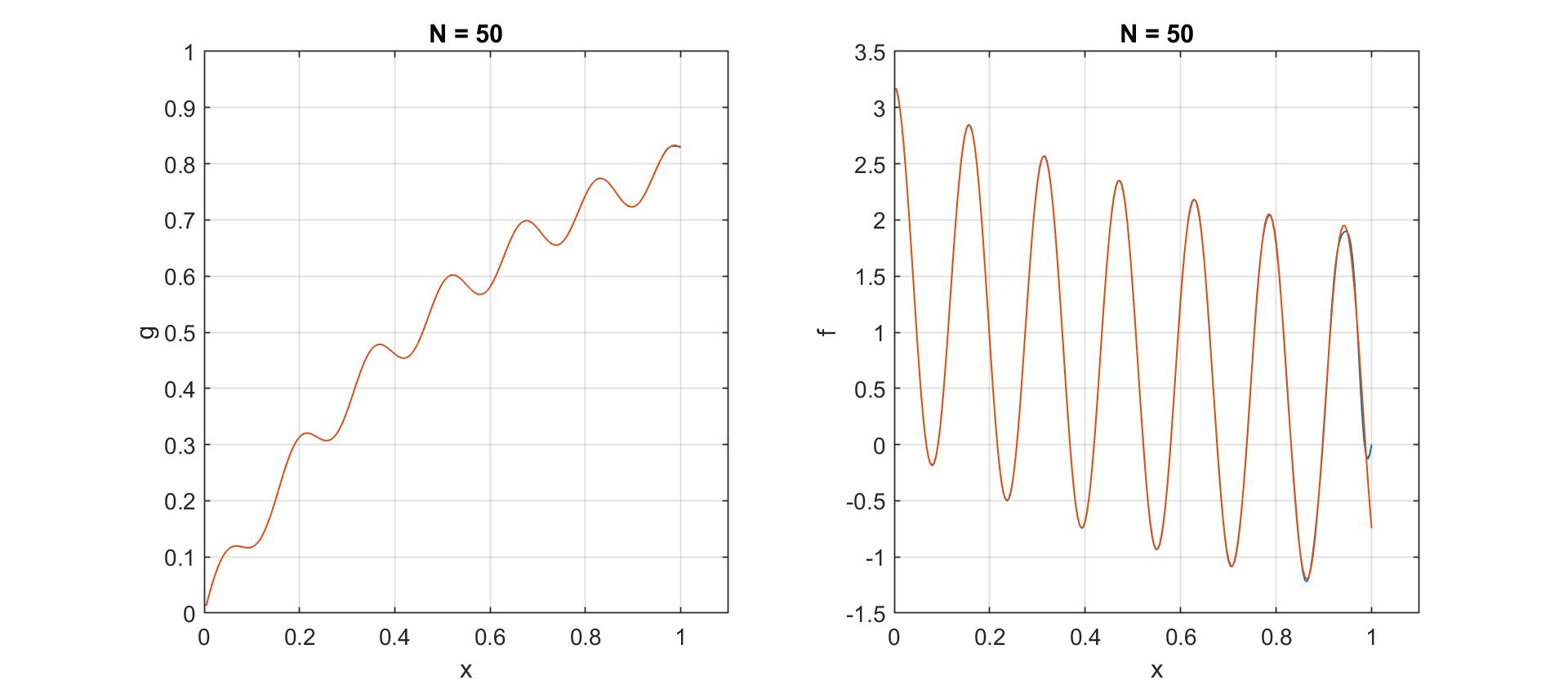}
   \caption{Plots of exact (red) $g$  \req{gcraig2}  (left) and $f$ \req{fcraig2}  (right) and respective approximations (blue) using sums of first $50$ terms of singular functions in \req{svdint}.  }
  \label{fig:CraigTrigUnpertbeta04N50}
\end{figure}

In our second example, just as in Sec.~\ref{sec:rustreg}, we perturb each data point $g(x_j)$ in the preceding example by obtaining a pseudorandom sample from a standard normal distribution, using the Matlab function $\mathbf{randn}$, multiplying it by a fixed standard deviation of $s = 0.05$, and then adding this number to $g(x_j)$; the result is depicted on the left-hand side in Fig.~\ref{fig:gpertfpert}. 
We then separate signal from noise in the data, using $\tau = 3$ as the threshold, as described in Sec.~\ref{sec:newreg}. The results, which are summarized in 
Fig.~\ref{fig:CraigTrigavectortau3}, are very similar to those in Fig.~\ref{fig:aprojMidPt}. Once again, the signal portion of $\mathbf{a}$ consists of the components  
$1,\,2,\,3,\,13,\,24,\,\mbox{and}\, 192$, but now the sum of squared residuals is $262.2$, which is slightly larger than the $258.1$ that was found using the midpoint rule in Sec.~\ref{sec:rustreg}. We again argue that components $a_{24}\ \mbox{and}\ {a}_{192}$ should be included in the noise, and we get the results that are summarized in Fig.~\ref{fig:Craig_Trig_gf_4term} and \ref{fig:Craig_Trig_stats_4term}. Even though the results are similar to those that we have already obtained using the midpoint rule, we have obtained much more information about both the source and the data functions. We now also have simple, closed-form, low-degree trigonometric polynomial approximations for both functions.

\begin{figure}[tbp] 
  \centering
  \includegraphics[width=5in,height=2.89in,keepaspectratio]{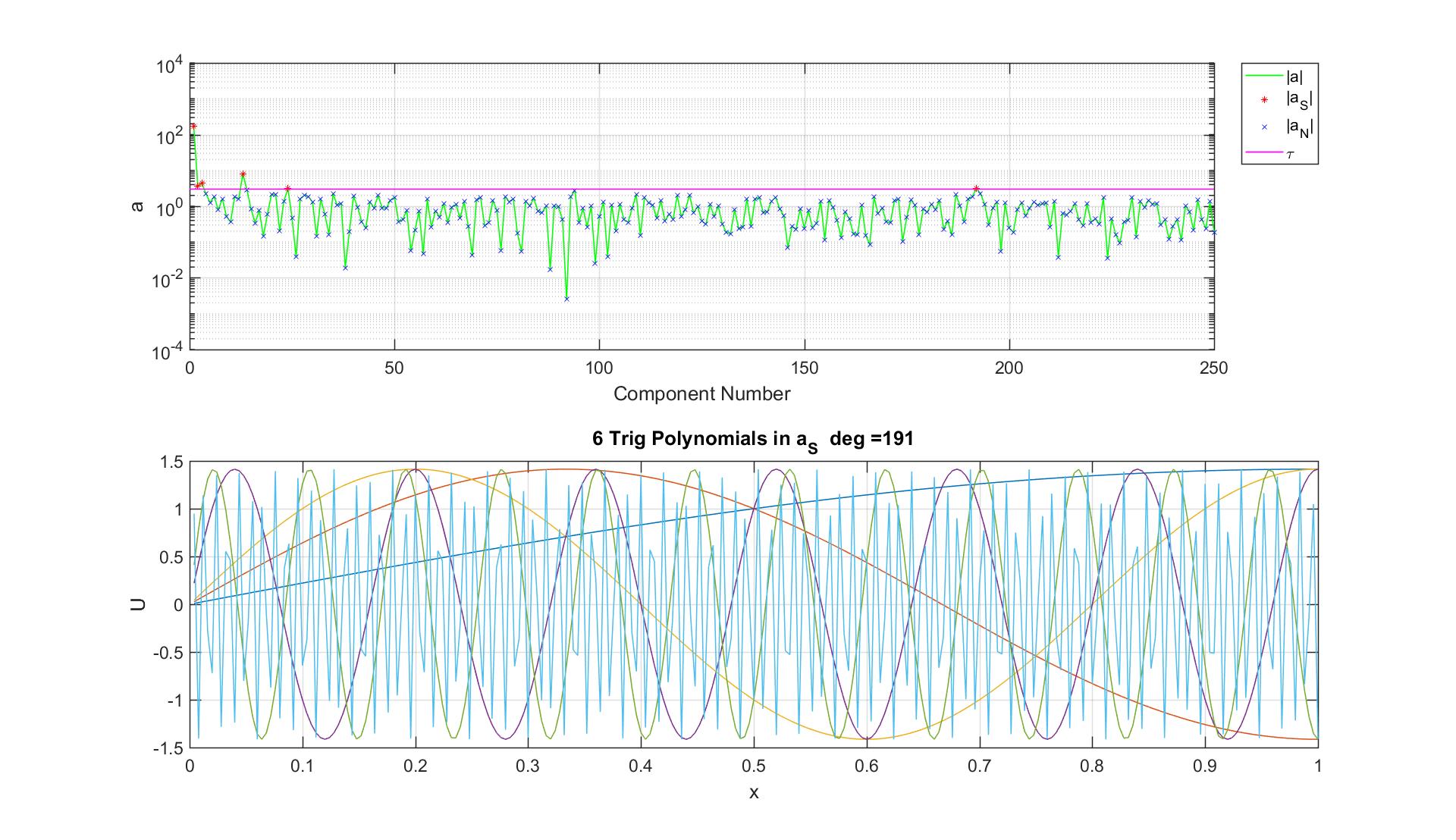}
  \caption{Separation of signal from noise in rotated data vector $\mathbf{a}$ (upper); columns of $\mathbf{U}$ that correspond to signal components in $\mathbf{a}$ (lower); note similarity to Fig.~\ref{fig:aprojMidPt} in Sec.~\ref{sec:rustreg}.}
  \label{fig:CraigTrigavectortau3}
\end{figure}

\begin{figure}[tbp] 
  \centering
  \includegraphics[width=5in,height=2.89in,keepaspectratio]{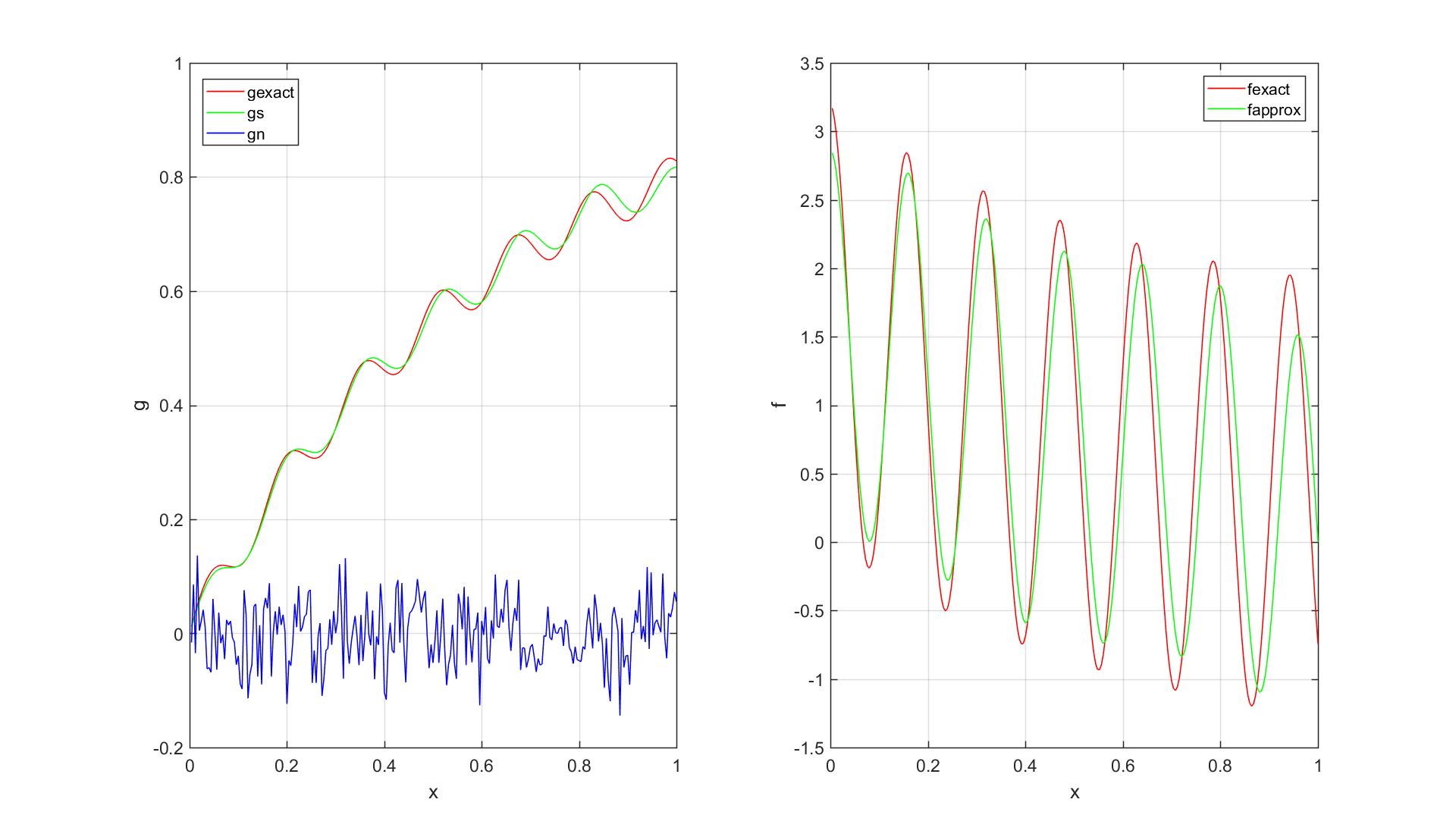}
  \caption{Truncation of signal from noise, and final approximations to data and source functions, with $a_{24}\ \mbox{and}\ {a}_{192}$ included in noise.}
  \label{fig:Craig_Trig_gf_4term}
\end{figure}

\begin{figure}[tbp] 
  \centering
  \includegraphics[width=5in,height=2.89in,keepaspectratio]{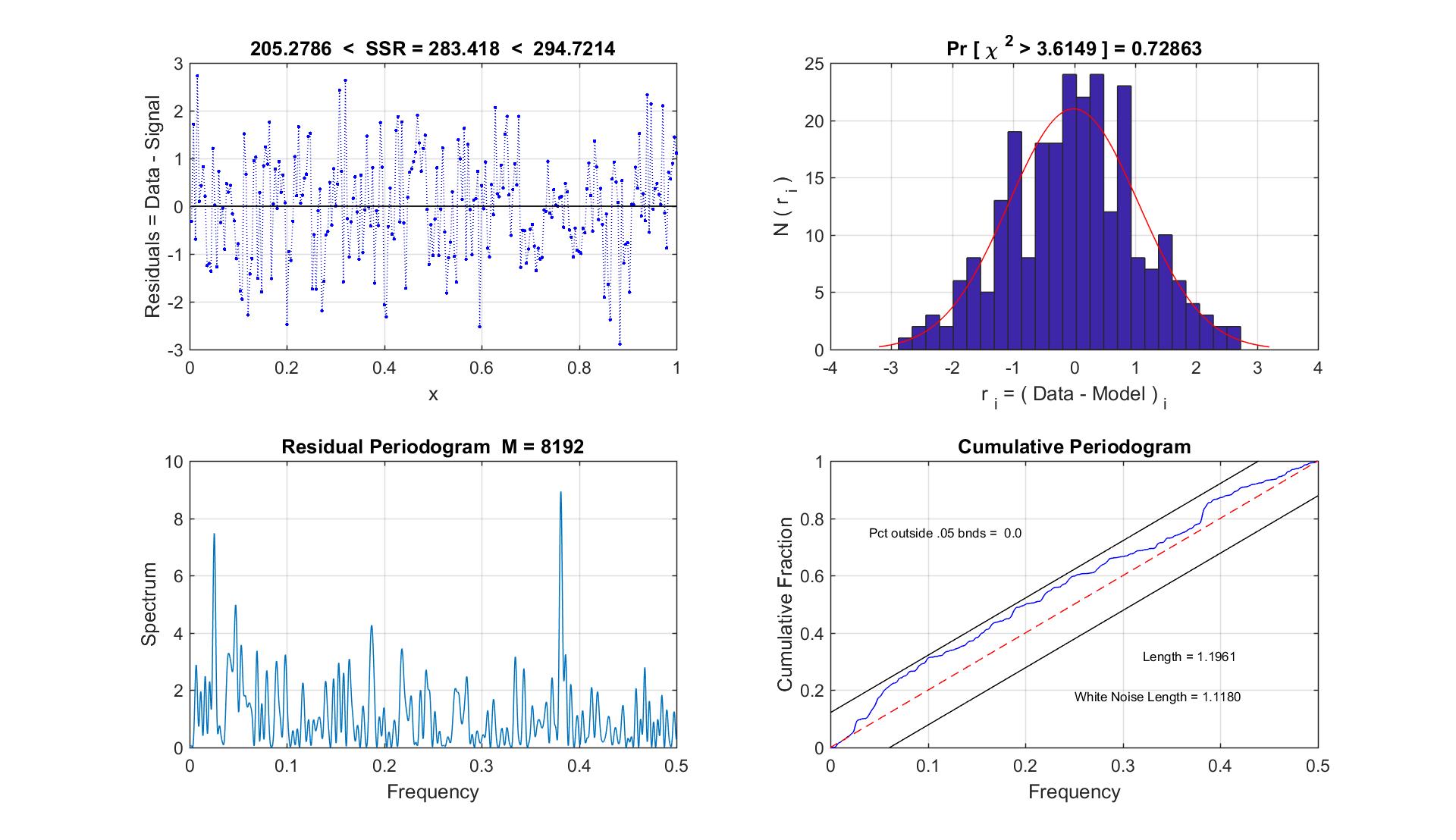}
  \caption{Diagnostics associated with noise in Fig.~\ref{fig:Craig_Trig_gf_4term} (left).}
  \label{fig:Craig_Trig_stats_4term}
\end{figure}

\subsection{Differentiating a Cubic Function}\label{subsec:cubic}

For our next example, we will estimate the derivative of the cubic function 
\be g(x) =  \frac{1}{2}\, [1+(2x-1)^3],\hspace{5ex}\frac{dg}{dx}(x) = 3\,(2x-1)^2,\hspace{5ex}0\le x\le 1,\labeq{cubic}\ee
on the same $250$-point grid as before. Before we perturb $g$, we once again check on using the trigonometric SVD expansions for the integration operator, which are given in Sec.~\ref{subsec:integration}, for $N=15,\,25,$ and $50$. The differences between the exact and approximate solutions are plotted in Fig.~\ref{fig:CubicfgTrig_unpert_error}. Once again, there is nonuniform convergence in the approximations of $g$, and a Gibbs phenomenon at $x=1$ in the approximations of $f$. Plots of the approximations of $g$ and $f$ for $N=50$ are given in Fig.~\ref{fig:CubicfgTrig_unpert_50terms}. Away from the right boundary, the approximation to $f$ provides a reasonable approximation of the source function.

\begin{figure}[tbp] 
  \centering
  \includegraphics[width=5in,height=2.35in,keepaspectratio]{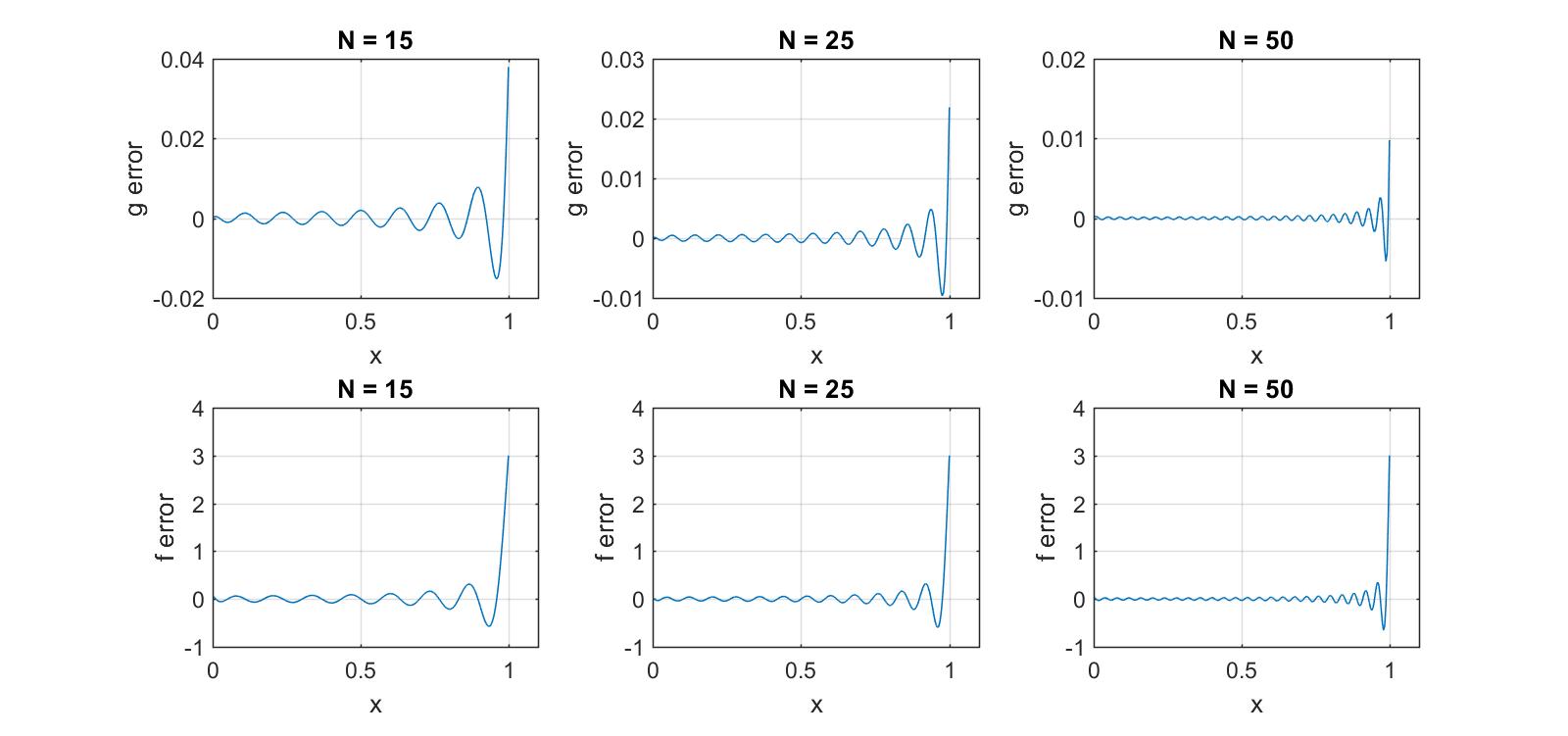}
  \caption{Differences between exact $g$ (top), and exact $f$ (bottom), for example with cubic $g$ \req{cubic}, and respective approximations using sums from $1$ to $N$ of trigonometric SVD expansions in Sec.~\ref{subsec:integration}.}
  \label{fig:CubicfgTrig_unpert_error}
\end{figure}

If we perturb the data points $g_j$ as above, using the same noise data and $\tau = 3$ as before, we get that the signal components of $\mathbf{a}$ consist of  
$1,\, 2,\,3,\,4,\,5,\,6,\,7,\,24,$ and $192$; see Fig.~\ref{fig:CubicfgTrig_a_truncation}. Moving $a_{24}$ and $a_{192}$ into the noise, as before, we get approximations of $g$ and $f$ that consist of the first seven terms in the respective trigonometric polynomial approximations of $G$ \req{gminterp} and $\hat{f}$ \req{hatfgen}. The results are given in Fig.~\ref{fig:CubicTrig7termapprox}. Here, the fact that $v_j(1) = 0$ for each $j$ has a significant effect on the estimate of the source function near the right endpoint.

\begin{figure}[tbp] 
  \centering
  \includegraphics[width=5in,height=2.52in,keepaspectratio]{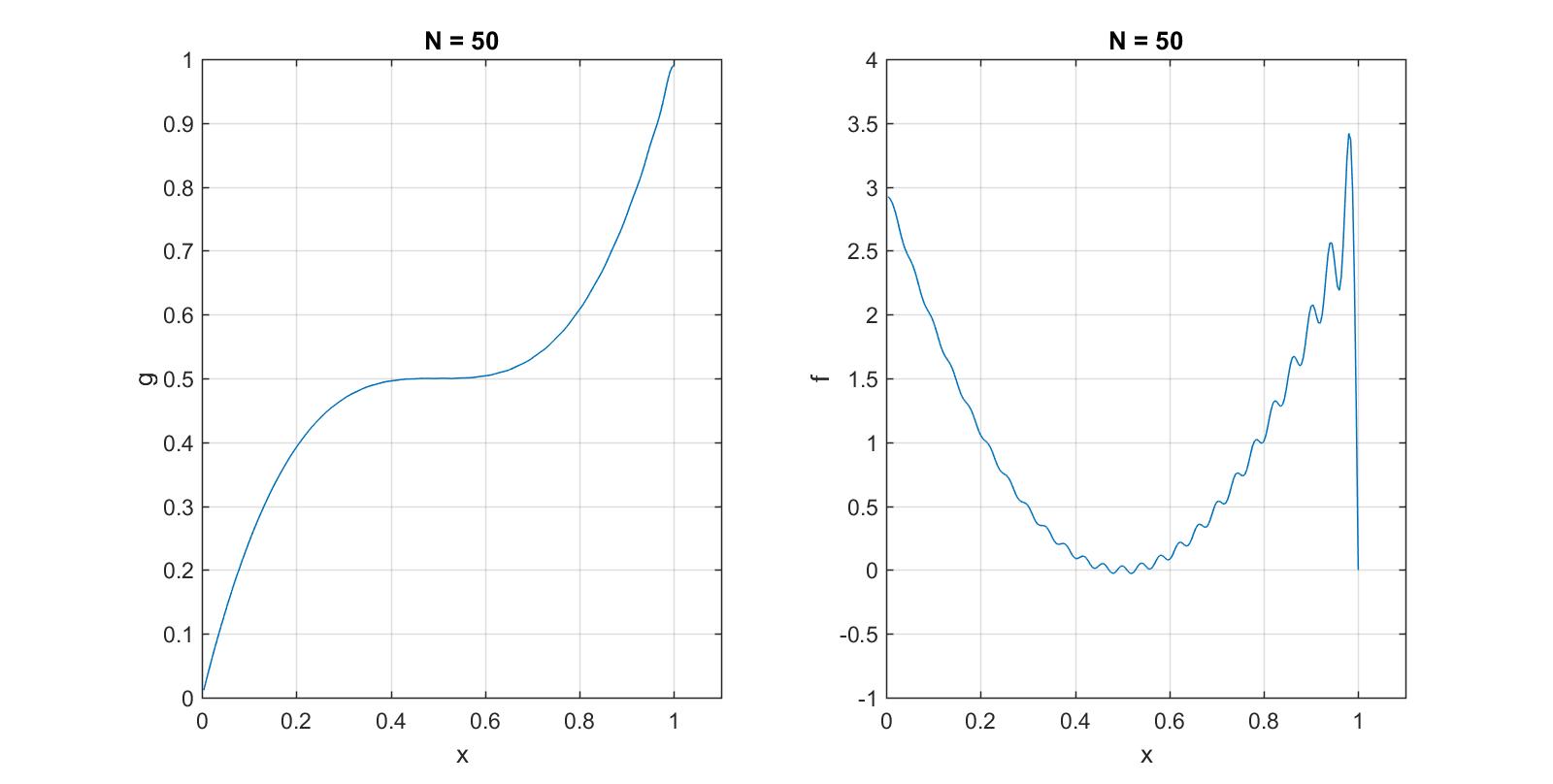}
 \caption{Trigonometric polynomial approximations to $g$ and $f$, with $N=50$, for unperturbed cubic example \req{cubic}.}
  \label{fig:CubicfgTrig_unpert_50terms}
\end{figure}

\begin{figure}[tbp] 
  \centering
  \includegraphics[width=5in,height=2.89in,keepaspectratio]{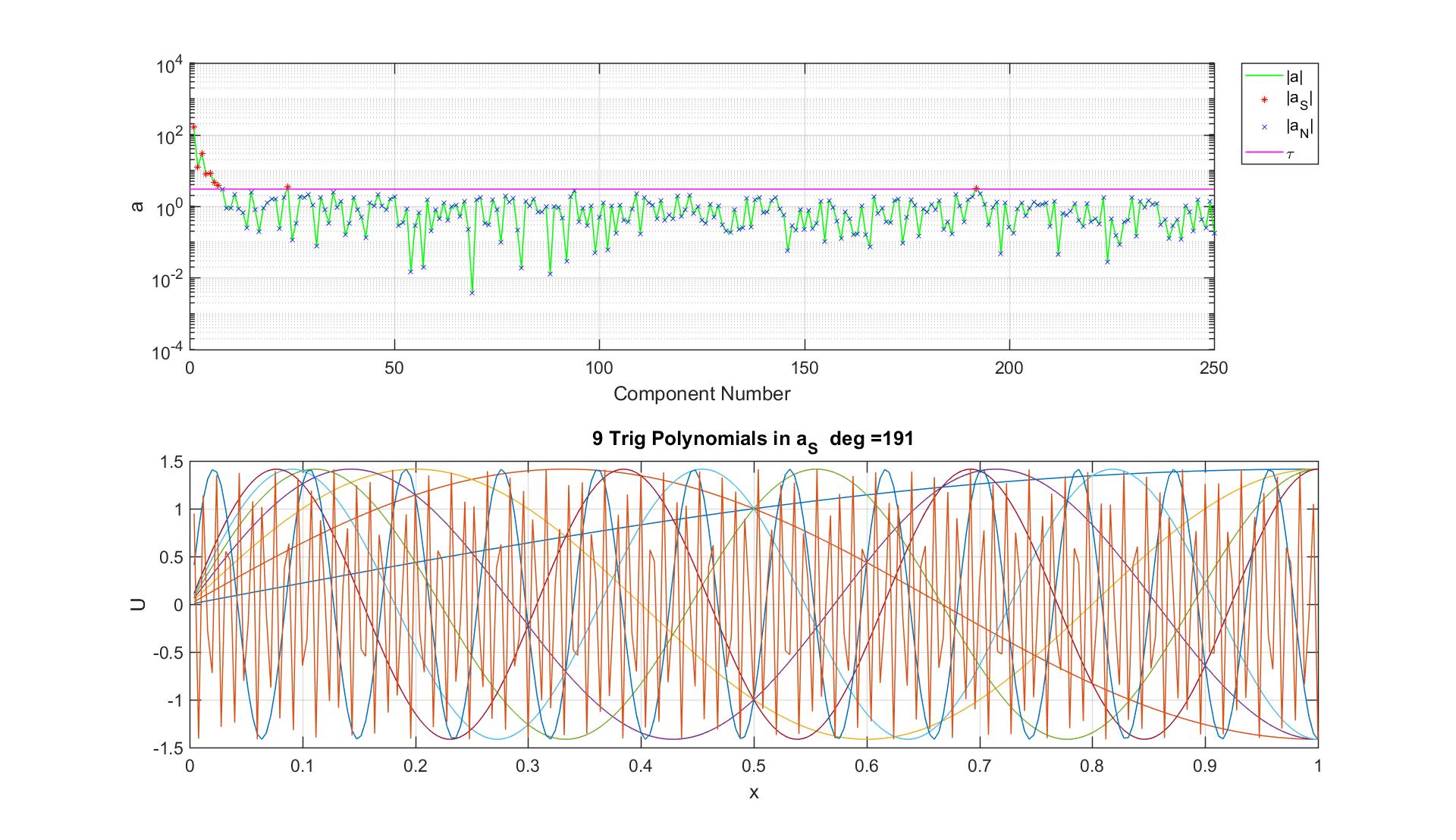}
   \caption{Truncation of signal from noise in $\mathbf{a}$ for cubic example \req{cubic}, perturbed by same noise vector as in example of Craig and Brown (upper); columns of $\mathbf{U}$ matrix that correspond to signal components in $\mathbf{a}$ (lower).}
  \label{fig:CubicfgTrig_a_truncation}
\end{figure}

\begin{figure}[tbp] 
  \centering
  \includegraphics[width=5in,height=2.89in,keepaspectratio]{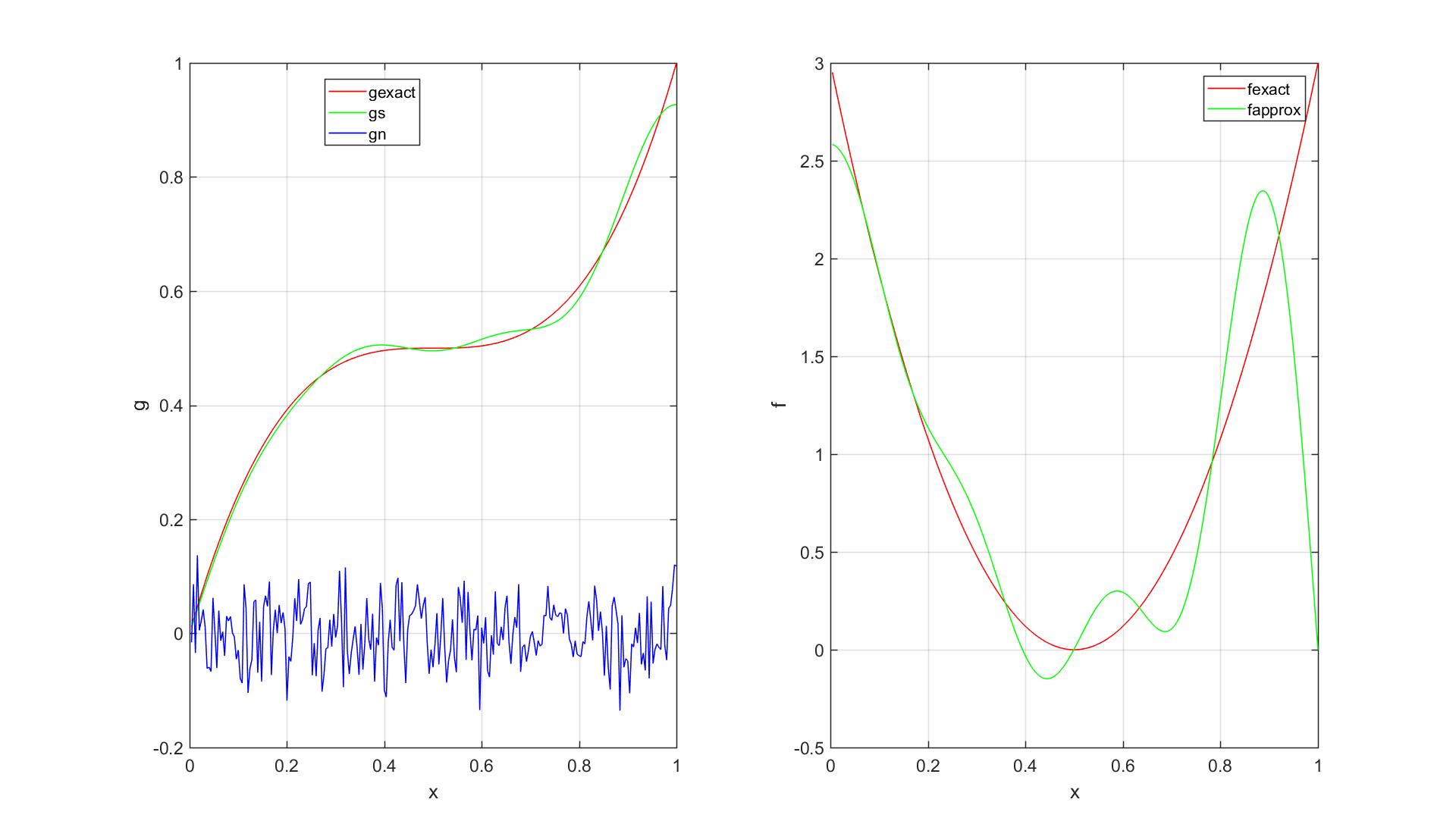}
 \caption{Seven term trigonometric polynomial approximations to $g$ and $f$ for cubic example \req{cubic}, with $a_{24}$ and $a_{192}$ included in noise.}
  \label{fig:CubicTrig7termapprox}
\end{figure}

\begin{figure}[tbp] 
  \centering
  \includegraphics[width=5in,height=2.89in,keepaspectratio]{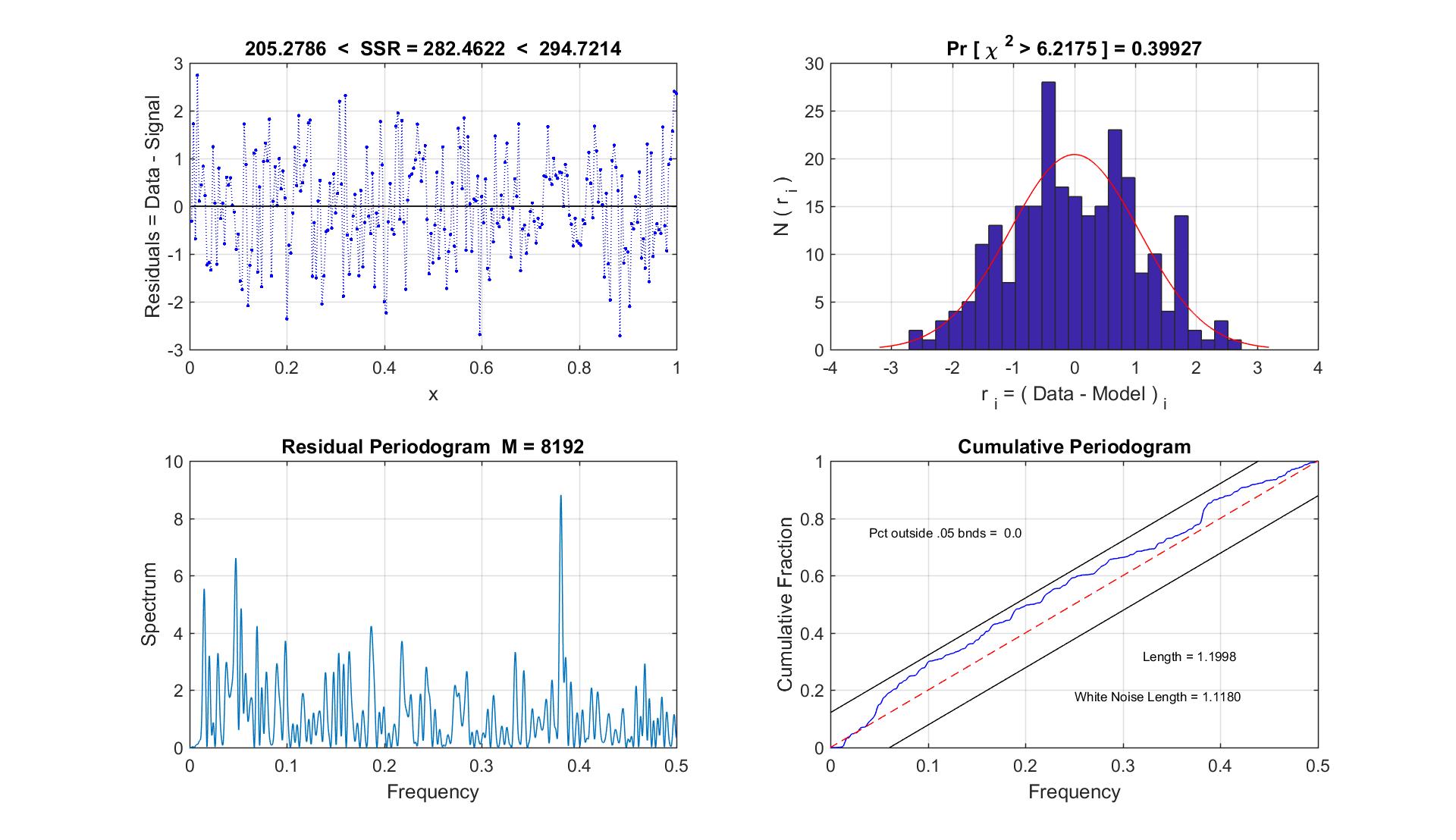}
   \caption{Residual diagnostics associated with noise vector in Fig.~\ref{fig:CubicTrig7termapprox} (left).}
  \label{fig:CubicTrig7termapprox_stats}
\end{figure}

It is interesting to regularize the perturbed data using a different set of basis functions, as has been discussed in Sec.~\ref{subsec:otherapps}. We transform Eq.~\req{cubic} to $[-1,1]$,
\be g(x) =  \frac{1}{2}\, (1+x^3),\hspace{5ex}\frac{dg}{dx}(x) = \frac{3}{2}\, x^2,\hspace{5ex}-1\le x\le 1,\labeq{cubic2}\ee
and instead of the sine polynomials on $[0,1]$ that are determined by the SVD of \req{integration}, we use Legendre polynomials, for the same noise  perturbation of the data $\mathbf{g}$ as above. Again using $\tau =3$ as a cutoff, we get the separation of signal from noise that is shown in Fig.~\ref{fig:CubicfgPleg_a_truncation}. Here, $\mathbf{a}_S$ consists of components $1,\,2,\,4,$ and $36$. Once again, we make the decision to include the higher-frequency component $a_{36}$ in the noise. Since there are closed-form expressions for the derivatives of the Legendre polynomials, we can estimate $\hat{f}$ by simply differentiating the expansion for $g$,
\be G(x) = \xi_1 P_0(x)+\xi_2 P_1(x)+\xi_4 P_3(x),\labeq{Plegcubicexpan}\ee
i.e., by computing
\be \hat{f}(x) = \xi_1 \frac{dP_0}{dx}(x) +\xi_2 \frac{dP_1}{dx}(x)+\xi_4 \frac{dP_3}{dx}(x).\labeq{DPlegcubicexpan}\ee
This gives the results in Fig.~\ref{fig:CubicfgPleg_approx3term} and \ref{fig:CubicfgPleg_approx3term_stats}. We emphasize two important points here. The first is that regularization of the noisy data function can be treated independently of the associated inverse problem of estimating the source function. The second point is that, once the data have been regularized, there may be better methods for estimating the source function than using the SVD expansion.

\begin{figure}[tbp] 
  \centering
  \includegraphics[width=5in,height=2.89in,keepaspectratio]{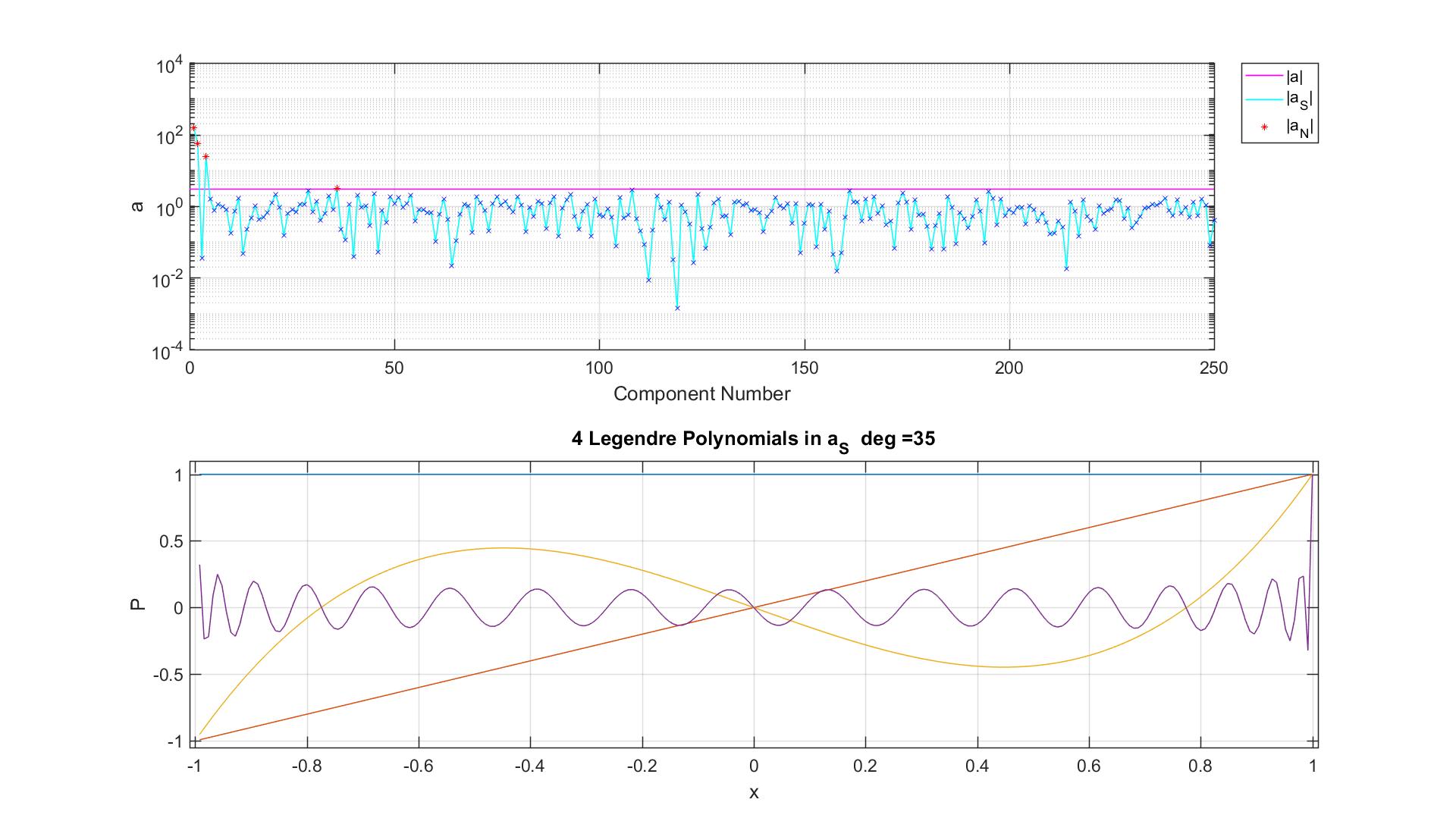}
  \caption{Separation of signal from noise in rotated data vector $\mathbf{a}$ for cubic example \req{cubic}, which has been regularized using Legendre polynomials. Plot of columns of $\mathbf{U}$ that correspond to signal components in $\mathbf{a}$.}
  \label{fig:CubicfgPleg_a_truncation}
\end{figure}

\begin{figure}[tbp] 
  \centering
  \includegraphics[width=5in,height=2.89in,keepaspectratio]{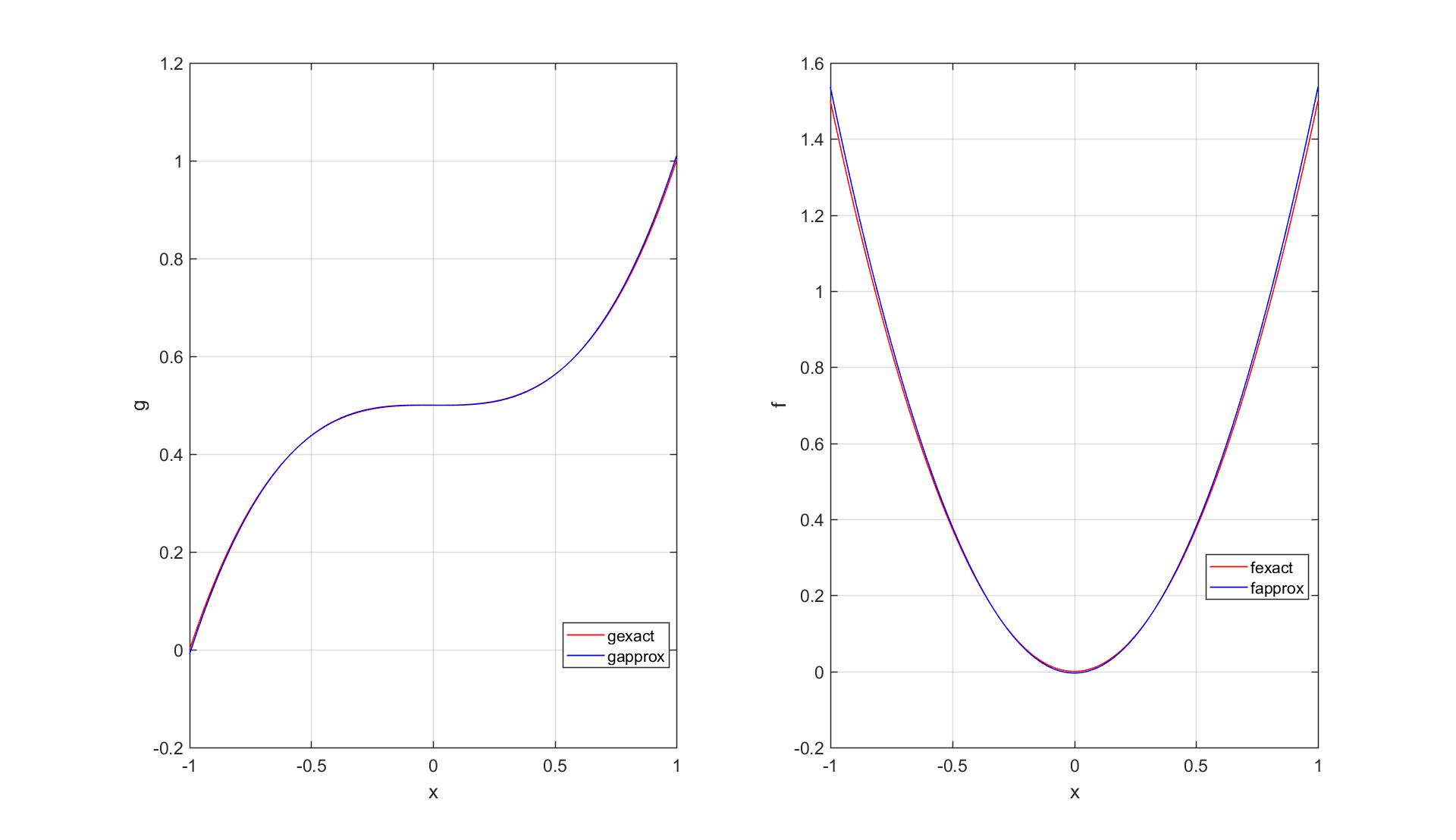}
  \caption{Final approximations for data and source using three-term Legendre polynomial approximation to $g$ and its first derivative approximation to $f$.}
  \label{fig:CubicfgPleg_approx3term}
\end{figure}

\begin{figure}[tbp] 
  \centering
  \includegraphics[width=5in,height=2.89in,keepaspectratio]{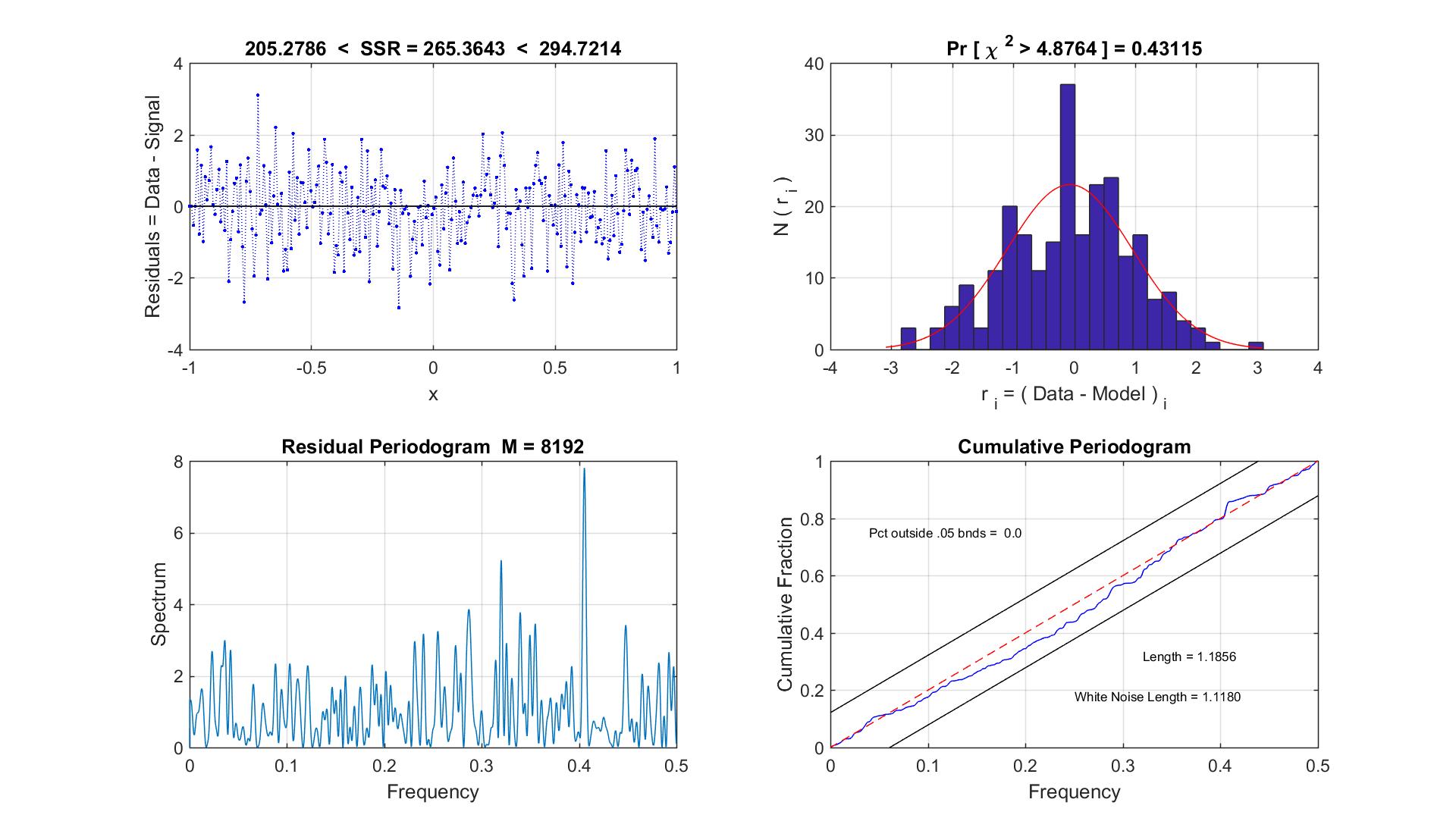}
  \caption{Diagnostics of the residual associated with three-term Legendre polynomial regularization of perturbed cubic function $g$ in \req{cubic}.}
  \label{fig:CubicfgPleg_approx3term_stats}
\end{figure}

\subsection{Fractional Differentiation Example}
For our next example, we consider an instance of Abel's equation \req{Abel} on $[-1,1]$, with $\mu = \nicefrac12$, which is the most frequent case in applications (\citep{craigbrown}, \citep{gorenflo+vessella}, \citep{golberg}, \citep{bracewell}), for which the exact solution is available \citep{Yousefi}, 
\bea g(x) &=& \frac{1}{\sqrt\pi}\frac{2}{105}\sqrt{\frac{x+1}{2}}\left[ 105-56\left(\frac{x+1}{2}\right)^2+48\left(\frac{x+1}{2}\right)^3\right],\labeq{Abelg}\\\vspace{4em}
f(x) &=&\frac{1}{\sqrt 2}\left[\left(\frac{x+1}{2}\right)^3 - \left(\frac{x+1}{2}\right)^2 +1\right].\labeq{Abelf}\eea
We note that $g(-1)=0$, and there is a singularity in $dg/dx$ at $x=-1$. Once again, we begin by checking the SVD expansions for fractional differentiation, as described in Sec.~\ref{subsec:fracint}. The approximations to $g$ and $f$ using the first $4$ terms are essentially indistinguishable; see Fig.~\ref{fig:Yousefi_Abel_fg_N4_20Oct17}.

\begin{figure}[tbp] 
  \centering
  \includegraphics[width=5in,height=2.59in,keepaspectratio]{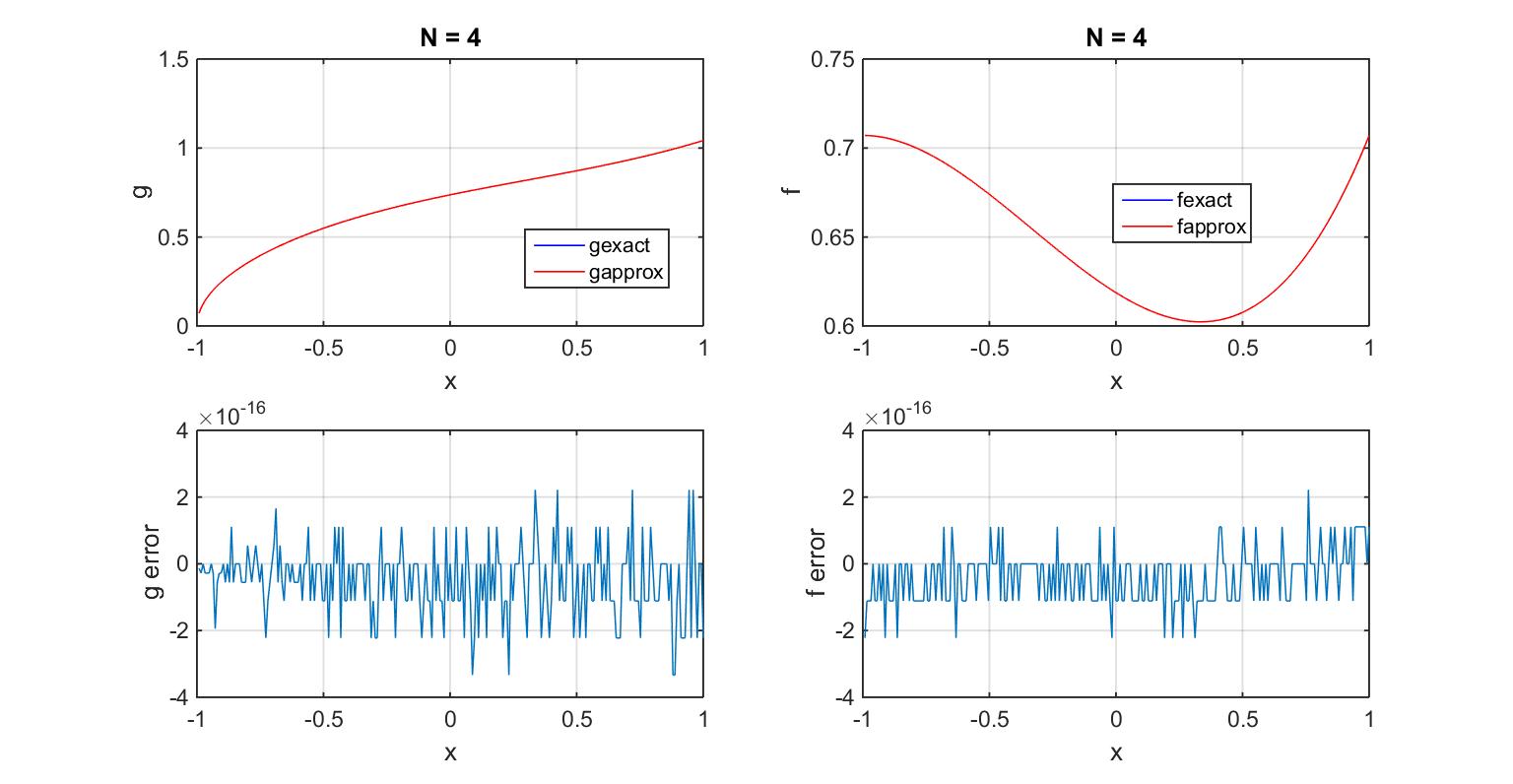}
  \caption{For Abel equation example \req{Abelg}-\req{Abelf}, exact $g$ and $f$ are essentially indistinguishable from approximations using first four terms of SVD expansions in Sec.~\ref{subsec:fracint}.}
  \label{fig:Yousefi_Abel_fg_N4_20Oct17}
\end{figure}

As our last example, we evaluate $g$ \req{Abelg} at $250$ equally-spaced points $-1<x_1<\cdots < x_{250} = 1$, and then we perturb these data at each point, as before, by computing a pseudorandom sample from a standard normal distribution, using the Matlab function $\mathbf{randn}$, multiplying this number by a fixed standard deviation of $s = 0.05$, and then adding the result to $g_j$. Proceeding to separate signal from noise as before, we encounter a difficulty with performing the $QR$ decomposition of the matrix $\mathbf{P}$ \req{JacobiC}, because its condition number blows up as the number of polynomials exceeds about $90$. The workaround we use, which could have been used on all of the earlier examples in which $g$ is perturbed by noise, takes advantage of {Diagnostic 4} in Sec.~\ref{sec:rustreg}, which is the assumption that high frequency components of $\mathbf{a}_S$ should be included in $\mathbf{a}_N$. Thus, instead of computing the full matrix $\mathbf{P}$, we compute the first $90$ columns of the matrix, determine the $QR$ decomposition of this smaller matrix, and then determine the signal components of $\mathbf{a}=\mathbf{Q}^T\mathbf{b}$ (see \req{wQTb}) using the resulting smaller $\mathbf{Q}$ matrix, again with the cutoff $\tau=3$. The result is summarized in Fig.~\ref{fig:Yousefi_aProjection20Oct17}. The signal components identified this way are $1,\,     2,\,     3,\,    87$. Using {Diagnostic 4} once more, we decide to include the higher-frequency component in $\mathbf{a}_N$. Using the remaining three components to determine the approximations to $g$ and $f$, we get the results in 
Fig.~\ref{fig:Yousefi_Abel_fg_approx20Oct17} and \ref{fig:Yousefi_Abel_fg_approx20Oct17_stats}. Once again, we have determined closed-form approximations, which are easy to compute, for the source and data functions. 

As in the case of integration, once the noisy data function has been regularized, there may be a more accurate method for estimating the source function. For example, since we have a closed-form expression for the regularized data function, this function could be evaluated on a denser mesh, that need not be equally spaced, and then a finite-difference estimate could be obtained for $f$ (e.g., \citet{linz}, \citet{brunner}). However, in this case, the estimate would be discrete, and not closed-form.  
  
\begin{figure}[tbp] 
  \centering
  \includegraphics[width=5in,height=2.89in,keepaspectratio]{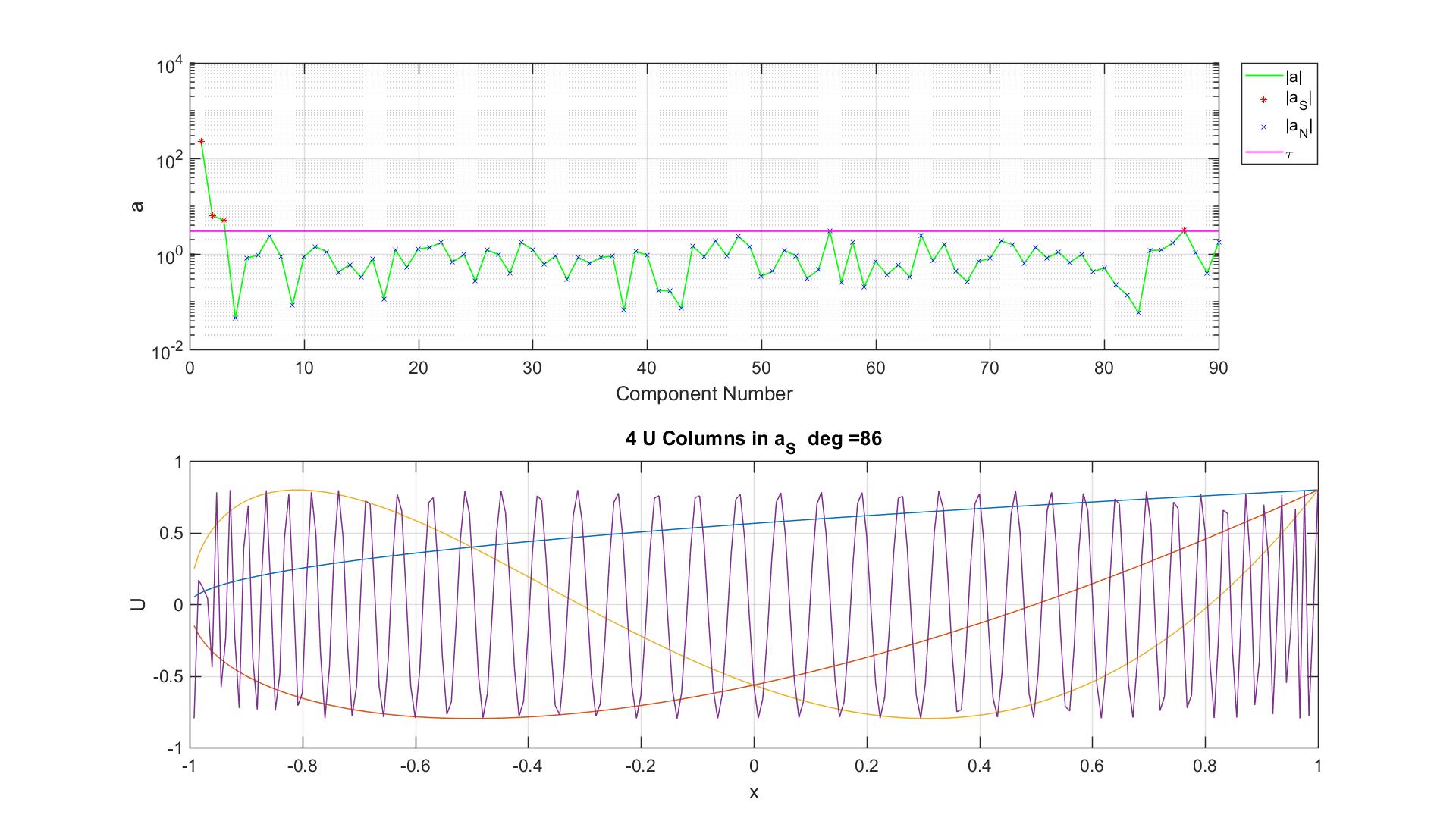}
 \caption{Separation of signal from noise in rotated data vector $\mathbf{a}$ for \req{Abelg}-\req{Abelf}, regularized using SVD expansion in Sec.~\ref{subsec:fracint} (upper); columns of $\mathbf{U}$ that correspond to signal components in $\mathbf{a}$ (lower).}
  \label{fig:Yousefi_aProjection20Oct17}
\end{figure}

\begin{figure}[tbp] 
  \centering
  \includegraphics[width=5in,height=2.89in,keepaspectratio]{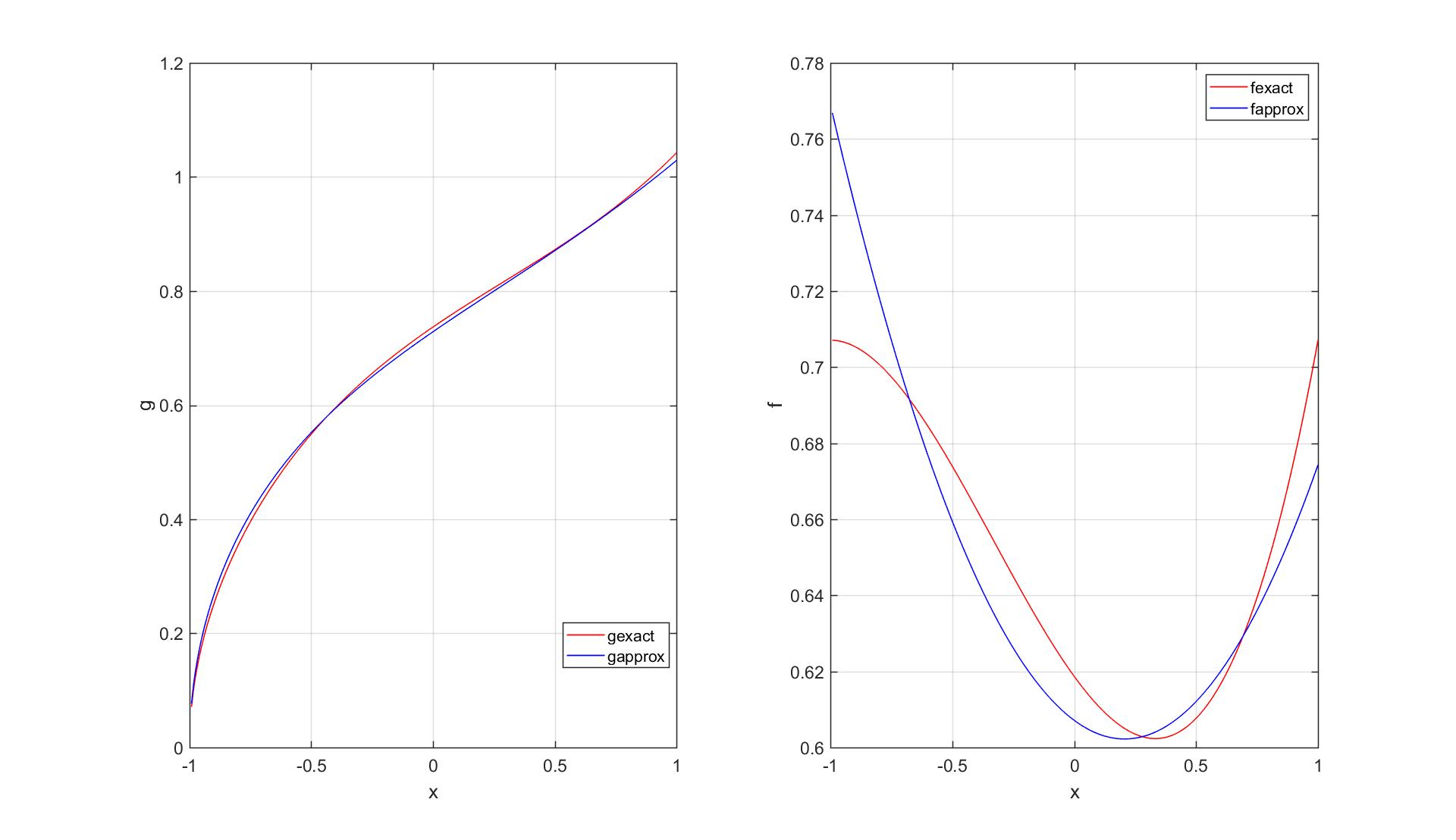}
 \caption{Three-term closed-form approximations of regularized data function \req{Abelg} and estimate of source function \req{Abelf} using SVD expansions in Sec.~\ref{subsec:fracint}.}
  \label{fig:Yousefi_Abel_fg_approx20Oct17}
\end{figure}

\begin{figure}[tbp] 
  \centering
  \includegraphics[width=5in,height=2.89in,keepaspectratio]{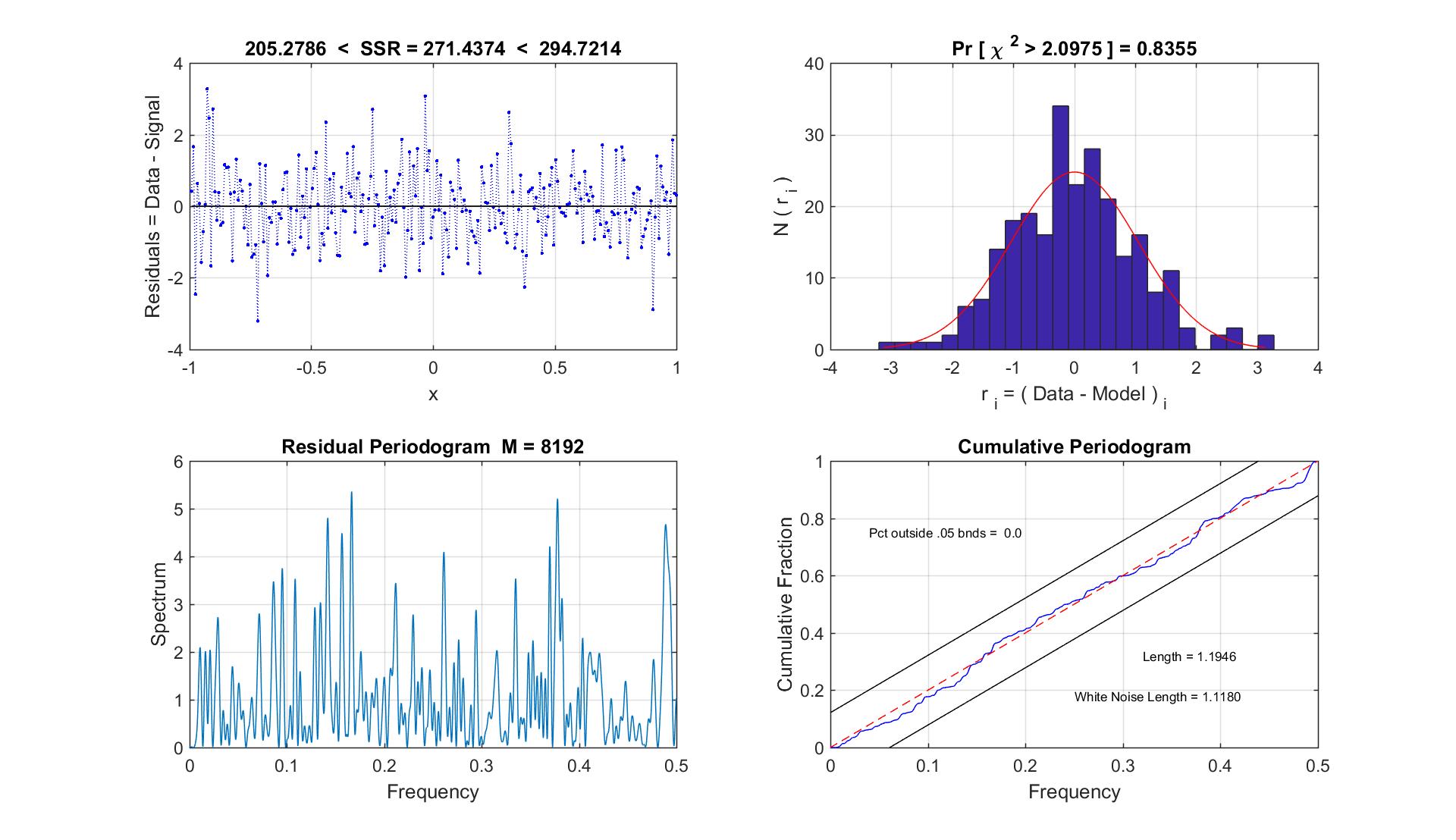}
 \caption{Diagnostics of residual associated with three-term regularization of perturbed data function \req{Abelg}.}
  \label{fig:Yousefi_Abel_fg_approx20Oct17_stats}
\end{figure}

\section{Concluding Remarks}\label{sec:conclusions}

We have presented a new method for regularizing a function that has been defined by a noisy set of measurement data. Our new approach extends a statistical time series method for separating signal from noise in the data that was introduced by Rust, by taking advantage of known closed-form singular value decompositions of the Hilbert space operators of integration and fractional integration. We have shown how to obtain a smooth, closed-form approximation of the data function, in the form of a linear combination of a small number of functions that are either trigonometric polynomials, Jacobi polynomials, or functions which are closely related to Jacobi polynomials, all of low degree. 
We have also shown how to obtain closed-form estimates of the derivative and fractional derivative of the data function, which are finite linear combinations of trigonometric or Legendre polynomials of low degree. 

\vspace{1ex}
We would like to acknowledge helpful discussions with Zydrunas Gimbutas and Howard Cohl, and we would like to thank Howard Cohl and Charles Hagwood for detailed comments on the manuscript.  

\vspace{1em}
{\it This paper is an official contribution of the National Institute of Standards and Technology and is not subject to copyright in the United States. Commercial products are identified in order to adequately specify certain procedures. In no case does such identification imply recommendation or endorsement by the National Institute of Standards and Technology, nor does it imply that the identified products are necessarily the best available for the purpose.}

\bibliographystyle{plainnat}
\bibliography{bertbib}

\begin{thebibliography}{24}
\providecommand{\natexlab}[1]{#1}
\providecommand{\url}[1]{\texttt{#1}}
\expandafter\ifx\csname urlstyle\endcsname\relax
  \providecommand{\doi}[1]{doi: #1}\else
  \providecommand{\doi}{doi: \begingroup \urlstyle{rm}\Url}\fi

\bibitem[gum(2008)]{gum}
Guide to the expression of uncertainty in measurement (gum), 2008.
\newblock URL \url{http://www.bipm.org/en/publications/guides/gum.html}.

\bibitem[Anderson(1958)]{anderson}
T.~W. Anderson.
\newblock \emph{An Introduction to Multivariate Statistical Analysis}.
\newblock John Wiley and Sons, New York, 1958.

\bibitem[Anderssen and de~Hoog(1990)]{golberg}
R.~S. Anderssen and F.~R. de~Hoog.
\newblock Abel integral equations.
\newblock In Michael~A. Goldberg, editor, \emph{Numerical Solution of Integral
  Equations}, chapter~8. Plenum Press, New York, New York, 1990.

\bibitem[Bracewell(2000)]{bracewell}
R.~N. Bracewell.
\newblock \emph{The Fourier Transform and its Applications}.
\newblock McGraw-Hill, Boston, third edition, 2000.

\bibitem[Brunner and van~der Houwen(1986)]{brunner}
H.~Brunner and P.~J. van~der Houwen.
\newblock \emph{The Numerical Solution of Volterra Equations}.
\newblock North-Holland, Amsterdam, New York, 1986.

\bibitem[Burkardt(2017)]{fsu+codes}
J.~Burkardt.
\newblock Matlab source codes, 2017.
\newblock URL \url{https://people.sc.fsu.edu/~jburkardt/m_src/m_src.html}.

\bibitem[Craig and Brown(1986)]{craigbrown}
I.~J.~D. Craig and J.~C. Brown.
\newblock \emph{Inverse Problems in Astronomy}.
\newblock Adam Hilger, Ltd, Bristol and Boston, 1986.

\bibitem[Cuer(1987)]{cuer}
M.~Cuer.
\newblock An illustrative problem on {A}bel's equation.
\newblock In P.~C. Sabatier, editor, \emph{Basic Methods of Tomography and
  Inverse Problems}, Malvern Physics Series, 4, pages 643--667. Adam Hilger,
  Bristol and Philadelphia, 1987.

\bibitem[{\relax DLMF}()]{NIST:DLMF}
{\relax DLMF}.
\newblock {\it NIST Digital Library of Mathematical Functions}.
\newblock http://dlmf.nist.gov/, Release 1.0.16 of 2017-09-18.
\newblock URL \url{http://dlmf.nist.gov/}.
\newblock F.~W.~J. Olver, A.~B. {Olde Daalhuis}, D.~W. Lozier, B.~I. Schneider,
  R.~F. Boisvert, C.~W. Clark, B.~R. Miller and B.~V. Saunders, eds.

\bibitem[Fuller(1976)]{fuller}
W.~A. Fuller.
\newblock \emph{Introduction to Statistical Time Series}.
\newblock John Wiley and Sons, New York, 1976.

\bibitem[Gander(1980)]{gander}
W.~Gander.
\newblock Algorithms for the qr-decomposition.
\newblock Technical Report 80-02, Eidgenoessische Technische Hochschule,
  Zurich, Switzerland, April 1980.

\bibitem[Gorenflo and Tuan(1995)]{GorenfloTuan1995}
R.~Gorenflo and K.~V. Tuan.
\newblock Singular value decomposition of fractional integration operators in
  ${L}_2$-spaces with weights.
\newblock \emph{Journal of Inverse and Ill-posed Problems}, 3:\penalty0 1--10,
  1995.

\bibitem[Gorenflo and Vessella(1991)]{gorenflo+vessella}
R.~Gorenflo and S.~Vessella.
\newblock \emph{Abel Integral Equations: Analysis and Applications}.
\newblock Lecture Notes in Mathematics (Book 1461). Springer, Berlin,
  Heidelberg, New York, London, 1991.

\bibitem[Hansen et~al.(2006)Hansen, Kilmer, and Kjeldsen]{ha+ki+kj}
P.~C. Hansen, M.~E. Kilmer, and R.~H. Kjeldsen.
\newblock Exploiting residual information in the parameter choice for discrete
  ill-posed problems.
\newblock \emph{BIT Numerical Mathematics}, 46:\penalty0 41--59, 2006.

\bibitem[Hogg and Craig(1965)]{hogg+craig}
R.~V. Hogg and A.~T. Craig.
\newblock \emph{Introduction to Mathematical Statistics}.
\newblock Macmillan, New York, second edition, 1965.

\bibitem[Kress(2014)]{kress}
R.~Kress.
\newblock \emph{Linear Integral Equations}.
\newblock Applied Mathematical Sciences. Springer, Berlin, Heidelberg, New
  York, London, third edition, 2014.

\bibitem[Linz(1985)]{linz}
P.~Linz.
\newblock \emph{Analytical and Numerical Methods for Volterra Equations}.
\newblock SIAM, Philadelphia, Pennsylvania, USA, 1985.

\bibitem[Morozov(1966)]{morozov1966}
V.~A. Morozov.
\newblock On the solution of functional equations by the method of
  regularization.
\newblock \emph{Doklady Mathematics}, 7\penalty0 (3):\penalty0 414 --417, 1966.
\newblock ISSN 1064-5624.

\bibitem[Rust(1998)]{BertIR}
B.~W. Rust.
\newblock Truncating the singular value decomposition for ill-posed problems.
\newblock Technical Report NISTIR 6131, National Istitute of Standards and
  Technology, Gaithersburg, MD, July 1998.

\bibitem[Rust(2000)]{rust2000}
B.~W. Rust.
\newblock Parameter selection for constrained solutions to ill-posed problems.
\newblock \emph{Computing Science and Statistics}, 32:\penalty0 333--347, 2000.

\bibitem[Rust and O'Leary(2008)]{rust+oleary}
B.~W. Rust and D.~P. O'Leary.
\newblock Residual periodograms for choosing regularization parameters for
  ill-posed problems.
\newblock \emph{Inverse Problems}, 24\penalty0 (3):\penalty0 1--30, 2008.
\newblock \doi{10.1088/0266-5611/24/3/034005}.

\bibitem[Snedecor and Cochran(1989)]{sn+co}
G.~W. Snedecor and W.~G. Cochran.
\newblock \emph{Statistical Methods}.
\newblock Iowa State University Press, eighth edition, 1989.

\bibitem[Wing(1991)]{wing}
G.~M. Wing.
\newblock \emph{A Primer on Integral Equations of the First Kind}.
\newblock SIAM, Philadelphia, Pennsylvania, USA, 1991.

\bibitem[Yousefi(2006)]{Yousefi}
S.~A. Yousefi.
\newblock Numerical solution of {A}bel's integral equation by using {L}egendre
  wavelets.
\newblock \emph{Applied Mathematics and Computation}, 175:\penalty0 574--580,
  2006.

\end{thebibliography}
\listoffigures
\end{document}